\documentclass{amsart}
\usepackage{amsfonts, amsmath, amssymb}
\usepackage{myarrows}
\usepackage[matrix, arrow, curve]{xy}

\DeclareMathOperator{\im}{im}
\DeclareMathOperator{\coker}{coker}

\DeclareMathOperator{\coeq}{coeq}

\def\Str{\mbox{\it String}}
\def\R{{\mathbb R}}
\def\Z{{\mathbb Z}}
\def\cL{{\mathcal L}}
\def\cC{{\mathcal C}}

\def\sk{\mathrm{sk}}
\def\str{\mathfrak{str}}
\def\g{\mathfrak g}

\def\simpl{{\mathrm{spl}}}
\def\proof{\noindent {\it Proof.} }
\def\qed{\hfill$\square$}
\def\intinf{{\textstyle \int}}

\def\inttwo{\tau_{\le 2}\intinf}
\def\dontshow#1
{}

\def\lllllarrow{\hspace{.05cm}\mbox{\,
\put(0,-3.5){$\leftarrow$}
\put(0,-.5){$\leftarrow$}
\put(0,2.5){$\leftarrow$}
\put(0,5.5){$\leftarrow$}
\put(0,8.5){$\leftarrow$}
               \hspace{.45cm}}}

\newtheorem{theorem}{Theorem}[section]
\newtheorem{proposition}[theorem]{Proposition} 
\newtheorem{definition}[theorem]{Definition} 
 
\newtheorem{lemma}[theorem]{Lemma}
\newtheorem{corollary}[theorem]{Corollary}
\newtheorem{assumptions}[theorem]{Assumptions}
\newtheorem{convention}[theorem]{\it Convention}
\newtheorem{remark}[theorem]{\it Remark}
\newtheorem{example}[theorem]{\it Example}

\begin{document}

\title{Integrating $L_\infty$-algebras}
\author{Andr\'e Henriques}

\maketitle      

\section{Introduction}

\subsection{Homotopy Lie algebras}

$L_\infty$-algebras, or strongly homotopy Lie algebras were introduced by Drinfeld and Stasheff \cite{Sta92} 
(see also \cite{HH93}), as a model for ``Lie algebras that satisfy Jacobi up to all higher homotopies''.

An $L_\infty$-algebra is a graded vector space $L=
L_0\oplus L_1\oplus L_2\oplus \cdots$ 
equipped with brackets \dontshow{brk}
\begin{equation}\label{brk}
[\,]:L\to L, \quad
[\,,]:\Lambda^2 L\to L, \quad
[\,,,]:\Lambda^3 L\to L,\quad 
\cdots 
\end{equation}
of degrees $-1$, $0$, $1$, $2\ldots$, where
the exterior powers are interpreted in the graded sense.
The various axioms satisfied by these brackets can be summarized as follows:

Let $L^\vee$ be the graded vector space with $L_{n-1}^*=Hom(L_{n-1},\R)$ in degree $n$, and let
$C^*(L):=Sym(L^\vee)$ be its symmetric algebra, again interpreted in the graded sense \dontshow{cocc}
\begin{equation}\label{cocc}
\begin{split}
C^*(L)=\R\oplus \Big[L_0^*\Big]&
\oplus \Big[\Lambda^2 L_0^*\oplus L_1^*\Big] \oplus \Big[\Lambda^3 L_0^*\oplus (L_0^*\!\otimes\! L_1^*)\oplus L_2^*\Big]\\
&\oplus \Big[\Lambda^4 L_0^*\oplus (\Lambda^2 L_0^*\!
\otimes\! L_1^*) \oplus Sym^2 L_1^*\oplus\ldots\Big]\oplus\ldots
\end{split}
\end{equation}
The transpose of the brackets (\ref{brk}) can be assembled into
a degree one map $L^\vee\to C^*(L)$, which 
extends uniquely to a derivation $\delta:C^*(L)\to C^*(L)$.
The equation $\delta^2=0$ then contains all the axioms that the brackets satisfy.

If the $k$-ary brackets are zero for all $k>2$, we recover the usual notion of differential graded Lie algebra.
If $L$ is concentrated in degrees $< n$, we get the notion of Lie $n$-algebra,
also called $n$-term $L_\infty$-algebra.
The $k$-ary brackets are then zero for all $k>n+1$.

The case of Lie 2-algebras has been studied by Baez and Crans \cite{BC04}.
A Lie 2-algebra consists of two vector spaces $L_0$ and $L_1$, and three brackets $[\,]$, $[\,,]$, $[\,,,]$ 
acting on $L=L_0\oplus L_1$, of degrees $-1,0,1$ respectively. 
A complete list of the axioms is given in \cite[Lemma 4.3.3]{BC04}.

Our main example, introduced by Baez and Crans \cite{BC04}, is the string Lie 2-algebra $\str=\str(\g):=\g\oplus \R$ 
associated to a simple Lie algebra $\g$ of compact type.
Its only non-zero brackets are
$[X,Y]\in\g$ for $X,Y\in\g$ and $[X,Y,Z]:=\langle[X,Y],Z\rangle\in\R$ for $X,Y,Z\in\g$.

\subsection{The string group and its Lie algebra}

The string group $\Str(n)$ is classically defined as the 3-connected cover of $\mbox{\it Spin}(n)$.
Extending the above definition, we shall call $\Str=\Str_G$ the 3-connected cover of any compact simple simply connected Lie group $G$.

Stolz and Teichner \cite{ST04}, \cite{Sto96} have a few different models for $\Str$,
one of which, inspired by work of Anthony Wassermann, 
is an extension of $G$ by the group of projective unitary operators in a type {\it III} von Neumann algebra.
Another model, due to Brylinski and McLaughlin \cite{BMcL}, consists of a $U(1)$-gerbe with connection over $G$.
More recently, John Baez et al \cite{BCSS05} came up with a model of $\Str$ in their quest for a Lie 2-group 
integrating the Lie 2-algebra $\str$.
Their model is a strict Lie 2-group, constructed using the loop group of $G$ and its basic central extension by $S^1$.

Strict 2-groups, i.e. group objects in the category of groupoids, are equivalent to crossed modules \cite{BL04}, \cite{For02}.
A crossed module consists of two groups $H_0$, $H_1$,
a homomorphism $\partial:H_1\to H_0$, and an action $H_0\acts H_1$ satisfying
\[
\partial(a\cdot b)=a\partial(b) a^{-1}\qquad\text{and}\qquad \partial(b)\cdot b' = b b' b^{-1},
\quad a\in H_0,\,\,\, b,b'\in H_1.
\] 
In the language of crossed modules, the Lie 2-group constructed in \cite{BCSS05} is given by
\(
\Str_0:=Path_*(G)$ and $\Str_1:=\widetilde{\Omega G},
\)
where $Path_*(G)$ denotes the based path space of $G$, and $\widetilde{\Omega G}$ the universal central extension of $\Omega G$.
The homomorphism $\partial$ is obtained by composing the projection $\widetilde{\Omega G}\to\Omega G$ with
the inclusion $\Omega G\hookrightarrow Path_*(G)$, while the action $Path_*(G)\acts\widetilde{\Omega G}$ is given by lifting
the pointwise conjugation action $Path_*(G)\acts\Omega G$.

It should be noted that the authors of \cite{BCSS05} did not provide a recipe for integrating Lie 2-algebras.
They simply constructed $\Str$, and observed that the corresponding (infinite dimensional) strict Lie 2-algebra
is equivalent to $\str$.
The present paper fills this gap by providing the desired integrating procedure.

To compare our results with those of \cite{BCSS05}, one replaces the crossed module $\Str$ by its simplicial nerve
\begin{equation*}
N\Str=\left[*\llarrow Path_*(G)\lllarrow \widetilde{Map_*(\partial \Delta^2, G)}\llllarrow \widetilde{Map_*(sk_1 \Delta^3, G)}\cdots\right].
\end{equation*}
That simplicial manifold agrees with our Lie 2-group (\ref{TheBstr}) integrating $\str$.
Here, the tilde over $Map_*(sk_1 \Delta^m, G)$ indicates the total space of an 
$(S^1)^{\left(\substack{\scriptstyle m\\ \scriptstyle 2}\right)}$-principal bundle.

\subsection{Homotopy Lie groups}

Let $\Delta$ be the category whose objects are the sets $[m]:=\{0,\ldots,m\}$, $m\ge 0$, and whose morphisms are the non-decreasing maps.
It is generated by $d^i:[m-1]\to[m]$ and $s^i:[m+1]\to[m]$ given by 
\begin{equation}\label{fdg}
\begin{matrix}
d^i(j)=j\phantom{+1}&\text{if}& j<i&&
s^i(j)=j\phantom{+1}&\text{if}& j\le i\\
d^i(j)=j+1&\text{if}& j\ge i&&
s^i(j)=j-1&\text{if}& j>i.
\end{matrix}
\end{equation}
\begin{definition}
A simplicial manifold $X$ is a 
contravariant functor from $\Delta$ to the category of manifolds.
We write $X_m$ for $X([m])$, and call $d_i:=X(d^i):X_m\to X_{m-1}$ and $s_i:=X(s^i):X_m\to X_{m+1}$
the face and degeneracy maps respectively.
We represent a simplicial manifold $X$ by a diagram
\[
X_0\llarrow X_1 \lllarrow X_2\llllarrow X_3\cdots
\]
in which the arrows stand for the face maps.
A simplicial manifold is called reduced if $X_0=pt$.
\end{definition}

Kan simplicial sets \cite[Section 1.3]{May92},
an important special class of simplicial sets,
are the simplicial sets $X$ such that any map $\Lambda[m,j]\to X$ extends to a map $\Delta[m]\to X$. 
Here $\Delta[m]$ denotes the simplicial $m$-simplex, 
and its $j$th horn $\Lambda[m,j]\subset\Delta[m]$ is the union of all its facets containing the $j$th vertex.
Similarly to the case of simplicial sets, we shall restrict our class of simplicial manifolds by imposing a version of the Kan condition.

\begin{definition}\label{hst} \dontshow{hst}
A simplicial manifold $X$ satisfies the Kan condition if the restriction maps \dontshow{dav}
\begin{equation}\label{dav}
X_m=Hom(\Delta[m],X) \to Hom(\Lambda[m,j],X)
\end{equation}
are surjective submersions for all $m$, $j$.
\end{definition}

\begin{remark}\rm
The above generalization of the Kan condition is different from that of Seymour \cite{Sey80}, in that he
requires the map (\ref{dav}) to have a {\em global} section, where as the condition we impose on it only implies the existence of local sections.
\end{remark}

It is well known that the category of groups can be embedded in the category of pointed spaces via the functor $G\mapsto K(G,1)$.
If $G$ is a Lie group, then $K(G,1)$ can be interpreted as a simplicial manifold
\begin{equation}\label{kg1}
K(G,1):=\big(*\llarrow G \lllarrow G^2 \llllarrow G^3\cdots\big),
\end{equation}
where the faces and degeneracies are given by the well known formulas
\[
\begin{split}
d_0(g_1,\ldots,g_m)&=(g_2,\ldots,g_m)\qquad
d_i(g_1,\ldots,g_m)=(g_1,\ldots,g_ig_{i+1},\ldots,g_m)\\
d_n(g_1,\ldots,g_m)&=(g_1,\ldots,g_{m-1})\quad
s_i(g_1,\ldots,g_m)=(g_1,\ldots,g_i,e,g_{i+1},\ldots,g_m).
\end{split}
\]
One can then characterize the reduced Kan simplicial manifolds (\ref{kg1}) coming from Lie groups
as those such that the maps (\ref{dav}) are diffeomorphisms for all $m>1$.
This justifies the following definition, which is a manifold analog of
Duskin's notion of $n$-group(oid) \cite{Dus79} (see also Glenn \cite{Gle82}). 

\begin{definition}\label{fk} \dontshow{fk}
A reduced simplicial manifold $X$ is called a Lie $n$-group if in addition to the Kan condition, 
the map (\ref{dav}) is a diffeomorphism for all $m>n$.
\end{definition}

The category of Lie 1-groups is equivalent to the usual category of Lie groups.
In an appendix, we relate our Lie 2-groups to those of Baez and Lauda \cite{BL04}:
there is a functor from their category to ours, which is fully faithful, but not essentially surjective.
We then provide a modification of their definition which fixes that lack of essential surjectivity.
For $n>2$, the notion of Lie $n$-group having not been defined in the literature, 
we don't need to further justify our definition.

\subsection{Integrating $L_\infty$-algebras}

The goal of this paper is to give a procedure that takes an $L_\infty$-algebra (a homotopy Lie algebra) 
and produces a reduced Kan simplicial manifold (a homotopy Lie group).
It assigns to $L$, the simplicial manifold $\intinf L$ given by \dontshow{in-pre}
\begin{equation}\label{in-pre}
\big(\intinf L\big)_m:=Hom_{DGA}\big(C^*(L),\Omega^*(\Delta^m)\big),
\end{equation}
where $C^*(L)$ denotes the Chevalley-Eilenberg complex (\ref{cocc}), and
$\Omega^*(\Delta^m)$ the differential graded algebra of de Rham differential forms on the $m$-simplex.
This is essentially the spatial realization construction of rational homotopy theory (see \cite{BG76}, \cite{Sul77}) applied to $C^*(L)$,
used by Hinich in the case of differential graded Lie algebras \cite{Hin97}, 
and by Getzler in the case of $L_\infty$-algebras \cite{Get04}.
If the $L_\infty$-algebra is nilpotent, Getzler then uses a gauge condition
to cut down (\ref{in-pre}) to an equivalent sub-simplicial space, which is finite dimensional in each degree.

In the absence of the nilpotence condition, we prefer not to modify (\ref{in-pre}) and rather to address 
the analysis inherent to the infinite dimensionality of these spaces.
More precisely, we concentrate on their manifold structure.
Provided we interpret $\Omega^*(\Delta^m)$ to mean
\[
\Omega^*(\Delta^m):=\big\{\alpha\,\big|\,\text{both $\alpha$ and $d\alpha$ are of class $C^r$}\big\},
\]
we prove that the spaces (\ref{in-pre}) are Banach manifolds, and that $\intinf L$ satisfies the Kan condition (Theorem \ref{tgo}).
We then set up a long exact sequence for computing the simplicial homotopy groups $\pi_n^\simpl(\intinf L)$ (Theorem \ref{tls}), and show that
the Lie algebra of $\pi_n^\simpl(\intinf L)$ is canonically isomorphic to $H_{n-1}(L)$.

If $L$ is a Lie $n$-algebra, we also provide a modification $\tau_{\le n}\intinf L$ of $\intinf L$
which satisfies the condition in Definition \ref{fk}.
It is not always a simplicial manifold, but 
we give a simple criterion which is equivalent to the spaces $(\tau_{\le n}\intinf L)_m$ being manifolds (Theorem \ref{TVF}).

\subsection{Acknowledgments}

I would like to thank Allen Knutson for telling me first about the String group, 
Jacob Lurie for explaining to me what $n$-groupoids are, 
Michael Hopkins for showing me simplicial sheaves, 
Ezra Getzler for a pleasant stay at Northwestern University and for a careful reading of this paper, 
Dmitry Roytenberg for conversations at an early stage of this work,
and Walter Paravicini for his help with Banach manifolds.

\section{Kan simplicial objects}\label{KSO}

In this sections, we fix a category $\cC$ equipped with a {\em Grothendieck pretopology}
and consider simplicial objects in that category.
We introduce the concepts of Kan simplicial objects,
and of Kan fibrations for maps between simplicial objects.
Finally, we introduce the simplicial homotopy groups of a Kan simplicial object,
and develop the theory of Postnikov towers.
Our formalism is closely related to Jardine's simplicial sheaves \cite{Jar86}.

We shall assume that $\cC$ has all coproducts, 
and use a slightly non-standard definition of Grothendieck pretopology.

\begin{definition}\label{GRO}
A Grothendieck pretopology on\, $\cC$ is a collection of morphisms, called covers,
subject the following three axioms:

Isomorphisms are covers.
The composition of two covers is a cover.
If $U\to X$ is a cover and $Y\to X$ is a morphism, then the pullback $Y\times_X U$ exists, and the natural morphism $Y\times_X U\to Y$ is a cover. 
\end{definition}

According to Grothendieck's original definition \cite[Expos\'e II]{SGA4I}, a cover should rather be a collection of morphisms $\{U_i\to X\}$.
Given a pretopology in our sense, one recovers a pretopology in Grothendieck's sense by declaring $\{U_i\to X\}$ to be a cover
if $\coprod_i U_i\to X$ is a cover.
See \cite[Section 0.3]{Joh77}, \cite[Chapter III]{MM94} for other discussions of Grothendieck pretopologies.

We make the following assumptions on our Grothendieck pretopology:

\begin{assumptions}\label{assu}
The category $\cC$ has a terminal object $*$, and for any object $X\in\cC$, the map $X\to *$ is a cover.

The category $\cC$ is Cauchy complete, which means that for any idempotent morphism $e:X\to X$, $e^2=e$,
there exists a subobject $i:A\to X$ and a retraction $r:X\to A$, $ri=1_A$, such that $ir=e$.

The retract of a cover is a cover. Namely, if we have a commutative diagram
\begin{equation}\label{retD}
\begin{matrix}
\xymatrix{
A\ar[r]\ar[d]^f&X\ar@/_8pt/[l]\ar[d]^g\\
B\ar[r]&Y\ar@/_8pt/[l]
}\end{matrix}
\end{equation}
where the composites $A\to X\to A$ and $B\to Y\to B$ are identities,
then $g$ being a cover implies that $f$ is also a cover.

Our pretopology is subcanonical, which means that all the representable functors $T\mapsto Hom(T,X)$ are sheaves.
\end{assumptions}

We'll sometimes want to talk about the limit of a diagram in $\cC$, before knowing that it actually exists.
For this purpose, we use the Yoneda functor 
\[
\begin{split}
\mathbf y\,:\,
\cC\,&\to \,\big\{\text{Sheaves on $\cC$}\big\}\\
X\,&\mapsto \big(T\mapsto hom(T,X)\big)
\end{split}
\]
to embed $\cC$ in the category of sheaves on $\cC$.
Using $\mathbf y$, a limit of objects of $\cC$ can always be interpreted as the limit of the corresponding representable sheaves.
The limit sheaf is then itself representable if and only if the original diagram had a limit in $\cC$.

Given a set $S$, and an object $X\in \cC$, we write $Hom(S,X)$ for the product $\prod_S X$, i.e. the object $X\times\ldots\times X$,
where the copies of $X$ are indexed by $S$.
If $S$ is now a simplicial set, and $X$ a simplicial object, we denote by $Hom(S,X)$ the obvious equalizer
\begin{equation}\label{SXSX}
Hom(S,X)\to \prod_{m\ge 0}Hom(S_m,X_m)\,\,\rrarrow \prod_{[m]\to[n]}Hom(S_n,X_m),
\end{equation}
which, as explained above, is a priori only a sheaf.
The sheaf $Hom(S,X)$ can also be described by giving its values on the various objects of $\cC$. 
Namely, for $T\in \cC$ we have
\begin{equation}\label{8bis}
\big(Hom(S,X)\big)(T):=hom(T\times S,X),
\end{equation}
where $T\times S$ denotes the simplicial object with $\coprod_{S_m}\!T$ in degree $m$.

Recall the Kan condition from Definition \ref{hst}.
The following generalizes it to other categories of simplicial objects.

\begin{definition}\label{dmk} \dontshow{dmk}
Let $\cC$ be a category equipped with a Grothendieck pretopology.
A map $f:X\to Y$ between reduced simplicial objects of $\cC$ is a Kan fibration if the obvious map from
$X_m=Hom(\Delta[m],X)$ to the object $Hom(\Lambda[m,j]\to\Delta[m],X\to Y)$ of commutative squares \dontshow{qsk}
\begin{equation}\label{qsk}
\begin{matrix}
\xymatrix@C=7mm@R=7mm{
\Lambda[m,j]\ar[r]\ar[d]&X\ar[d]^f\\
\Delta[m]\ar[r]&Y
}\end{matrix}
\end{equation}
is a cover.
A simplicial object $X$ is called Kan if the map $X\to *$ is a Kan fibration.
\end{definition}

The object of squares (\ref{qsk}) is a certain equalizer, similar to (\ref{SXSX}).
It is best described in the language of sheaves by saying that its value at $T$ is the set of commutative squares
\begin{equation}\label{Ttimes}
\begin{matrix}
\xymatrix@C=7mm@R=7mm{
T\times \Lambda[m,j]\ar[r]\ar[d]&X\ar[d]\\
T\times \Delta[m]\ar[r]&Y
}\end{matrix}
\end{equation}
with given vertical arrows.
The fact that the above sheaf is representable is not obvious from the definition. 
It relies on the following Lemma. 

Recall that a simplicial set $S$ is {\em collapsable} if it admits a filtration \dontshow{ss*}
\begin{equation}\label{ss*}
*=S_0\subset S_1 \subset\cdots\subset S_k=S
\end{equation}
such that each $S_i$ is obtained from the previous one by filling a horn, namely such that
$S_i$ can be written as $S_i=S_{i-1}\times_{\Lambda[n_i,j_i]}\Delta[n_i]$ for some map $\Lambda[n_i,j_i]\to S_{i-1}$.

\begin{lemma}\label{collapsable1}\dontshow{collapsable1}
Let $S\subset\Delta[n]$ be a collapsable simplicial set.
Let $f:X\to Y$ be a map between reduced simplicial objects that satisfies the conditions of Definition \ref{dmk} for all $m<n$
(this implicitly means that $Hom(\Lambda[m,j]\to\Delta[m],X\to Y)$ is representable).
Then the object of commutative squares $Hom(S\to\Delta[n],X\to Y)$ exists in\, $\cC$.
\end{lemma}

\proof 
Filter $S$ as in (\ref{ss*}).
We show by induction on $i$ that $Hom(S_i\to\Delta[n],X\to Y)$ exists in $\cC$.
We begin the induction by noting that $Hom(*\to\Delta[n],X\to Y)=Y_n$ is an object of $\cC$.

We now assume that $Hom(S_{i-1}\to\Delta[n],X\to Y)$ exists in $\cC$.
The bottom row of the pullback diagram
\[
\xymatrix{
Hom(S_i\to\Delta[n],X\to Y)\ar[r]\ar[d]\ar@{}[rd]|(.45){\line(0,1){5}\line(-1,0){5}}& Hom(S_{i-1}\to\Delta[n],X\to Y)\ar[d]\\
Hom(\Delta[n_i],X)\ar@{->}[r]& Hom(\Lambda[n_i,j_i]\to\Delta[n_i],X\to Y)
}
\]
is a cover by hypothesis. It follows that $Hom(S_i\to\Delta[n],X\to Y)$ exists in $\cC$.
\qed

Since the horns $\Lambda[n,j]$ are collapsable, we have:

\begin{corollary}\label{ikq}\dontshow{ikq}
Let $f:X\to Y$ be a map satisfying the conditions of Definition \ref{dmk}
for all $m<n$. Then the object $Hom(\Lambda[n,j]\to\Delta[n],X\to Y)$ of commutative squares (\ref{qsk}) exists in $\cC$.
\end{corollary}

We now collect a few easy lemmas to be used in the future.

\begin{lemma}\label{rwfj}\dontshow{rwfj} 
The category $\cC$ has all finite products; the product of two covers is a cover.
\end{lemma}
\proof
Since $X\to *$ is a cover, the pullback $X\times _* Y$ exists.
The product of two covers $X\to X'$, $Y\to Y'$ can be factorized as $X\times Y\to X\times Y'\to X'\times Y'$.
The first morphism is a cover because it is the pullback of $Y\to Y'$ along the projection $X\times Y'\to Y'$, and similarly for the second morphism.
The result follows since covers are closed under composition. 
\qed

\begin{lemma}\label{rmf}\dontshow{rmf} 
The product of two Kan fibrations of simplicial objects is again a Kan fibration.
\end{lemma}
\proof
This follows from the equality
$Hom(\Lambda[m,j]\to\Delta[m],X\times X'\to Y\times Y')=
Hom(\Lambda[m,j]\to\Delta[m],X\to Y)\times Hom(\Lambda[m,j]\to\Delta[m],X'\to Y')$, 
and the fact that the product of two covers is a cover.
\qed

\begin{lemma}\label{cmg}\dontshow{cmg}
The composite $X\to Y\to Z$ of two Kan fibrations of simplicial objects is again a Kan fibration.
\end{lemma}

\proof In the diagram below, the top left and the bottom horizontal arrows are covers.
The square being a pullback, it follows that the top composite is also a cover.
\[
\xymatrix{
X_m\ar[r]&Hom(\Lambda[m,j]\to\Delta[m],X\to Y)\ar[r]\ar[d]\ar@{}[rd]|(.45){\line(0,1){5}\line(-1,0){5}}& 
Hom(\Lambda[m,j]\to\Delta[m],X\to Z)\ar[d]\\
&Y_m\ar[r]& Hom(\Lambda[m,j]\to\Delta[m],Y\to Z)\\
}
\]
\qed

\begin{lemma}\label{lith}\dontshow{lith}
Let $A\rightarrow B$ be an acyclic cofibration between finite simplicial sets, and let $X\to Y$ be a Kan fibration of simplicial objects.
Assume that the object of commutative squares $Hom(A\to B,X\to Y)$ exists in\, $\cC$.
Then $Hom(B,X)$ exists in $\cC$, and the map
\begin{equation}\label{HBX}
Hom(B,X)\,\,\rightarrow \,\,Hom(A\to B,X\to Y)
\end{equation}
is a cover.
\end{lemma}

\proof
By \cite[Section IV.2]{GZ67} (see also \cite[Corollary 2.1.15]{Hov99}), there exists a simplicial set $S$ containing $B$ such that $B$ is a retract of $S$ via maps 
$j:B\to S$ and $r:S\to B$ fixing $A$, and such that $S$ admits a filtration
\(
A=S_0\subset S_1 \subset\cdots\subset S_k=S
\)
with $S_i=S_{i-1}\times_{\Lambda[n_i,j_i]}\Delta[n_i]$.
As in Lemma \ref{collapsable1}, we have pullback diagrams
\[
\xymatrix{
Hom(S_i\to B,X\to Y)\ar[r]\ar[d]\ar@{}[rd]|(.45){\line(0,1){5}\line(-1,0){5}}& Hom(S_{i-1}\to B,X\to Y)\ar[d]\\
Hom(\Delta[n_i],X)\ar@{->}[r]& Hom(\Lambda[n_i,j_i]\to\Delta[n_i],X\to Y)
}
\]
where the maps $S_i\to B$ are given by $r|_{S_i}$.
The bottom arrow is a cover, and therefore so is the top arrow.
Composing all the above maps, we conclude that $Hom(S\to B,X\to Y)$ exists in $\cC$ and that 
\begin{equation}\label{HSBX}
Hom(S\to B,X\to Y)\,\,\rightarrow\,\, Hom(A\to B,X\to Y)
\end{equation}
is a cover.
The sheaf $Hom(B,X)$ is a retract of $Hom(S\to B,X\to Y)$ via maps
\[
\Big(B\stackrel{f}{\to} X\Big)\mapsto
\left(\begin{matrix}S \stackrel{f\circ r}{\to} X\\ {\scriptstyle r}\downarrow\hspace{.7cm}\downarrow{\scriptstyle \pi}\\B \stackrel{\pi\circ f}{\to}Y\end{matrix}\right)
\quad\text{and}\quad
\left(\begin{matrix}S \stackrel{\phi}{\to} X\\ {\scriptstyle r}\downarrow\hspace{.7cm}\downarrow{\scriptstyle \pi}\\B \stackrel{\psi}{\to}Y\end{matrix}\right)
\mapsto
\Big(B\stackrel{\phi\circ j}{\to} X\Big),
\]
where for better readability, we have omitted all the ``\,$T\!\times\!\cdots$'' that should be there as in (\ref{8bis}) and (\ref{Ttimes}).
It follows that $Hom(B,X)$ is representable in $\cC$.
It is not hard to check that the above maps are compatible with the projections to $Hom(A\to B,X\to Y)$.
The morphism (\ref{HBX}) is therefore a retract of (\ref{HSBX}), and thus also a cover.
\qed

\begin{corollary}\label{cor-lith}\dontshow{cor-lith}
Let $A\to B$ be an acyclic cofibration between finite simplicial sets, and let $X$ be a Kan simplicial object.
Assume that $Hom(A,X)$ exists in $\cC$. Then $Hom(B,X)$ exists in $\cC$, and the map 
\[
Hom(B,X)\,\,\rightarrow \,\,Hom(A,X)
\]
is a cover.
In particular, if $B$ be a finite contractible simplicial set and $X$ a Kan simplicial object, 
then $Hom(B,X)$ exists in $\cC$. \qed
\end{corollary}

Without the assumption that the sheaf $Hom(A,X)$ is representable in $\cC$, we still have the following analog of Corollary \ref{cor-lith}.
\begin{lemma}\label{SURJ}
Let $A\to B$ be an acyclic cofibration between finite simplicial sets, and let $X$ be a Kan simplicial object.
Then the map 
\[
Hom(B,X)\,\,\rightarrow \,\,Hom(A,X)
\]
is a surjective map of sheaves.
\end{lemma}
\proof
In the statement of Corollary \ref{cor-lith}, replace $\cC$ by the category of sheaves on $\cC$, and the notion of cover by the notion of surjective map of sheaves.
\qed

\section{Simplicial homotopy groups}\label{Simplicial homotopy groups}

Recall \cite[Section I.4]{May92} that given a reduced Kan simplicial set $X$, its homotopy groups are given by
\begin{equation}\label{1dpi}
\pi_n(X):=\{x\in X_n\,|\,d_i(x)=*\text{ \rm for all } i\}/\sim,
\end{equation}
where $x\sim x'$ if there exists an element $y\in X_{n+1}$ such that $d_0(y)=x$, $d_1(y)=x'$, and $d_i(y)=*$ for all $i>1$.

More generally, let $S$ be any simplicial set whose geometric realization is homotopy equivalent to the $n$-sphere.
And let $SI$ be any simplicial set equipped with an injective map $S\vee S\to SI$ with the property that the two inclusion
$S\to SI$ are homotopy equivalences, and are homotopic to each other.
Then $\pi_n(X)$ is isomorphic to the coequalizer
\begin{equation}\label{2dpi}
\coeq\Big(Hom(SI,X)\,\rrarrow Hom(S,X)\Big).
\end{equation}
Definition (\ref{1dpi}) is the special case of (\ref{2dpi}) in which $S=\Delta[n]/\partial \Delta[n]$, and 
$SI=\Delta[n+1]\,\big/\big(\!\underset{i=2}{\overset{n}{\cup}}F^i\cup(F^0\cap F^1)\big)$,
where $F^i$ denotes the $i$-th facet of $\Delta[n+1]$.

Homotopy groups in the category of simplicial sheaves were introduced in \cite{Jar86}.
The following is essentially Jardine's definition.

\begin{definition}\label{S and SI}
Let $X$ be a reduced Kan simplicial object and let $S$, $SI$ be as above. 
The simplicial homotopy groups $\pi_i^{\simpl}(X)$ are the sheaves \dontshow{shg}
\begin{equation}\label{3dpi}
\pi_n^{\simpl}(X):=\coeq\Big(Hom(SI,X)\,\rrarrow Hom(S,X)\Big).
\end{equation}
The fact that this is independent of the choice of $S$ and $SI$ is the content of Lemma \ref{S'}.
\end{definition}

To make the above definition a little bit more concrete,
we write down what it means for the simplest case $S=\Delta[n]/\partial \Delta[n]$,
$SI=\Delta[n+1]\,\big/\!\big(\underset{i=2}{\overset{n}{\cup}}
F^i\cup(F^0\cap F^1)\big)$. In that case, we get
\[
\pi_n^{\simpl}(X)(T)=\left\{\begin{minipage}{2.4cm}
- a cover $U\to T$

- a morphism\\ $x:U\to X_n$\end{minipage}\,\left|
\,\begin{minipage}{6cm}
- $d_i\circ x=*$ for all $i$

- there exist $y:U\times_T U\to X_{n+1}$ 
such that 
$d_0\circ y=x\circ\mathrm{pr}_1$, $d_1\circ y=x\circ\mathrm{pr}_2$, and $d_i\circ y=*$ for all $i>1$.
\end{minipage}\right\}\right.\Big/\sim, 
\]
where $(U,x)\sim (U',x')$ if there exist a
covering $V$ of $U\times_X U'$ and a
map $y':V\to X_{n+1}$ such that 
$d_0\circ y'=x\circ\mathrm{pr}_1$, $d_1\circ y'=x'\circ\mathrm{pr}_2$, and $d_i\circ y'=*$ for all $i>1$.

To show that $\pi_n^\simpl(X)$ are groups in the category of sheaves, we
let $(U,x)$, $(V,y)$ be two representatives of elements of $\pi_n^{\simpl}(X)(T)$.
By lemma \ref{SURJ}, the map
\[
Hom(\Delta[n+1],X)\to Hom(\Lambda[n+1,1],X)
\]
is a surjective.
So we may pick a common refinement $U\leftarrow W\to V$, and a map $w:W\to X_{n+1}$
such that $d_0\circ w=x$, $d_2\circ w=y$ and $d_i\circ w=*$ for $i>2$.
The product is then represented by $(W,d_1\circ w)$.

\begin{example}\rm\label{kge}\dontshow{kge}
Let $G$ be a group object of $\cC$, and let $K(G,1)$ be as in (\ref{kg1}). 
Then $\pi_1^\simpl(K(G,1))$ is the sheaf represented by $G$, and $\pi_n^\simpl(K(G,1))=0$ for $n\ge 2$.
\end{example}

\begin{lemma}\label{S'}
Let $S$, $S'$ be simplicial sets weakly equivalent to the $n$-sphere, 
and let $S\vee S\to SI$, $S'\vee S'\to SI'$ be monomorphisms with the property that the two inclusions
$S\to SI$ are homotopic homotopy equivalences, and similarly for the two inclusions $S'\to SI'$.
Then the sheaves $\pi_n^{\simpl}(X):=\coeq(Hom(SI,X)\,\rrarrow Hom(S,X))$ and $\pi_n^{\simpl}(X)':=\coeq(Hom(SI',X)\,\rrarrow Hom(S',X))$ are canonically isomorphic.
\end{lemma}

\proof
Since $S$ and $S'$ are weakly equivalent, we may find a third simplicial set $S''$, and acyclic cofibrations $S\hookrightarrow S''\hookleftarrow S'$.
Similarly, we may find acyclic cofibrations $SI\hookrightarrow SI''\hookleftarrow SI'$ making the diagram
\[
\xymatrix{
S\vee S\ar[r]\ar[d]&
S''\vee S''\ar[d]&
S'\vee S'\ar[l]\ar[d]\\
SI\ar[r]&SI''&SI'\ar[l]
}
\]
commute. Letting $\pi_n^{\simpl}(X)'':=\coeq(Hom(SI'',X)\,\rrarrow Hom(S'',X))$,
the above diagram induces maps $\pi_n^{\simpl}(X)\to\pi_n^{\simpl}(X)''\leftarrow \pi_n^{\simpl}(X)'$.
We just show that the first map is an isomorphism, the argument being the same for the second one.

The left and middle vertical maps in
\[
\xymatrix{
Hom(SI'',X)\ar[d]^{(a)}\ar@<.7ex>[r]\ar@<-.7ex>[r]&Hom(S'',X)\ar[d]^{(b)}\ar[r]&\pi_n^{\simpl}(X)''\ar[d]^{(c)}\\
Hom(SI,X)\ar@<.7ex>[r]\ar@<-.7ex>[r]&Hom(S,X)\ar[r]&\pi_n^{\simpl}(X)
}
\]
are surjective by Lemma \ref{SURJ}. 
The surjectivity of $(b)$ implies the surjectivity of $(c)$,
and the surjectivity $(a)$ implies the injectivity of $(c)$.
\qed

All basic properties of homotopy groups of Kan simplicial sets go through to Kan simplicial objects.
For example, we have:

\begin{proposition}\label{sess}\dontshow{sess}
Let $X\to Y$ be a Kan fibration of Kan simplicial objects, and let $F$ be its fiber.
Then there is an associated long exact sequence 
\begin{equation}\label{LLLEEESSS}
\ldots\to\pi^\simpl_{n+1}(Y)\stackrel{\partial}{\to}\pi_n^\simpl(F)\to\pi_n^\simpl(X)\to\pi_n^\simpl(Y)\to\ldots
\end{equation}
\end{proposition}

\proof
At the level of points, the proof is standard and can be found in \cite[Theorem 7.6]{May92}.
At the level of sheaves, we just need to keep track of covers and refine them appropriately each time we use the Kan condition.

A more complete argument, including a precise construction of the connecting homomorphism $\partial$, can be found in \cite[Section II]{VOs77}.
\qed

We now introduce truncation functors that kill the simplicial homotopy groups outside of a given range.
Recall that for Kan simplicial sets, one has two truncation functors $\tau_{\le n}$ and $\tau_{<n}$ 
given by
\begin{equation*}
(\tau_{\le n}X)_m=X_m\big/\!\sim \qquad\text\qquad(\tau_{< n}X)_m=X_m\big/\!\approx
\end{equation*}
where $x\sim x'$ if the corresponding maps $\Delta[m]\to X$ are simplicially homotopic relatively to the $(n-1)$-skeleton 
of $\Delta[m]$, 
and $x\approx x'$ if they are equal when restricted to the $(n-1)$-skeleton 
of $\Delta[m]$.
The $i$th simplicial homotopy group of $\tau_{<n}X$ (of $\tau_{\le n}X$) is equal to $\pi_i^\simpl (X)$ for $i<n$ (respectively $i\le n$) and zero otherwise. 
In particular, $\tau_{<n+1}X \rightarrow\tau_{\le n}X$ is a weak equivalence.

The functor $\tau_{<n}$ is due to Moore and is by now classical \cite{May92} \cite{GJ99}.
The functor $\tau_{\le n}$ is apparently due to Duskin and is less frequent in the literature.
We refer the reader to \cite[Section 3.1]{Gle82}, \cite[Proposition 1.5]{Bek04}, \cite[Section 2]{Get04} for further discussions about it.

If one is willing to leave the category of simplicial objects of $\cC$ and allow simplicial sheaves, 
then one has analogs of the above truncations functors.

\begin{definition}\label{DEFTAU}
Let $X$ be a reduced Kan simplicial object.
Then the Postnikov pieces $\tau_{\le n}X$ and $\tau_{<n}X$ are the simplicial sheaves given by 
\begin{equation*}
(\tau_{\le n}X)_m=\coeq\Big(
Hom\big(P
,X\big)\,\rrarrow Hom\big(\Delta[m],X\big)\Big),
\end{equation*}
\begin{equation*}
(\tau_{< n}X)_m=\im\Big(Hom(\Delta[m],X)\rightarrow Hom(\sk_{n-1}\Delta[m],X)\Big),
\end{equation*}
where $P$ is the pushout simplicial set
\[
\xymatrix{
\Delta[1]\times\sk_{n-1}\Delta[m]\ar[r]\ar[d]\ar@{}[rd]|(.55){\line(0,-1){5}\line(1,0){5}}&\sk_{n-1}\Delta[m]\ar[d]\\
\Delta[1]\times\Delta[m]\ar[r]&\,P\,.
}
\]
\end{definition}

We then have a tower \dontshow{pKr}
\begin{equation}\label{pKr}
\ldots \to \tau_{<n+1}X \to \tau_{\le n}X \to \tau_{<n}X \to \ldots
\to \tau_{\le 1}X\to *
\end{equation}
whose inverse limit is $X$. 
As before, the $i$th simplicial homotopy group of $\tau_{<n}X$ (of $\tau_{\le n}X$) is equal to $\pi_i^\simpl (X)$ for $i<n$ (respectively $i\le n$), 
and zero otherwise. 

The truncation functor $\tau_{\le n}$ will later be used to construct Lie $n$-groups integrating Lie $n$-algebras.
We can already check that $\tau_{\le n}X$ satisfies one of the conditions in Definition \ref{fk}.

\begin{lemma}\label{check n-group}
Let $X$ be a Kan simplicial object, and let $\tau_{\le n}X$ be as above.
Then
\[
(\tau_{\le n}X)_m\to Hom\big(\Lambda[m,j],\tau_{\le n}X\big)
\]
is an isomorphism for all $m> n$.
\end{lemma}
\proof
We only prove the lemma in the case when the topos of sheaves on $\cC$ has enough points \cite[Section IX.11]{MM94}.
Since stalks of Kan simplicial objects are Kan simplicial sets, and since $\tau_{\le n}$ commutes with taking stalks,
it is enough to prove the lemma for the case of simplicial sets.

So let $X$ be a Kan simplicial set.
Let $S$ be the union of all $n$-faces of $\Delta[m]$ containing its $j$-th vertex.
Since $S$ is weakly homotopy equivalent to $\Delta[m]$ as simplicial set under $\sk_{n-1}\Delta[m]$,
and since $X$ is fibrant, the restriction map
\[
\begin{split}
(\tau_{\le n}X)_m=\big\{\text{maps $\Delta[m]\to X$ modulo homotopy fixing $\sk_{n-1}\Delta[m]$}\big\}&\\
\to\big\{\text{maps $S\to X$ modulo homotopy fixing $\sk_{n-1}\Delta[m]$}\big\}\hspace{.2cm}&
\end{split}
\]
is a bijection.
An element of $(\tau_{\le n}X)_m$ can thus be described as a map $f:\sk_{n-1}\Delta[m]\to X$,
along with a collection of homotopy classes of nullhomotopies, one for each restriction $f|_{\partial F}$ to an $n$-face $F$ of $S$. 

A map $\Lambda[m,j]\to \tau_{\le n}X$ is therefore the same thing as a map $\sk_{n-1}\Lambda[m,j]\to X$,
along with a collection of homotopy classes of nullhomotopies, one for each restriction $f|_{\partial F}$ to an $n$-face $F$ of $\Lambda[m,j]$ containing the
$j$th vertex. 
The result follows since $\sk_{n-1}\Lambda[m,j]=\sk_{n-1}\Delta[m]$ and since 
every $n$-face $F$ of $\Delta[m]$ containing its $j$th vertex is also contained in $\Lambda[m,j]$.
\qed

%

\begin{remark}\rm
For the standard argument that shows how to reduce the proof of Lemma \ref{check n-group} to the case of simplicial sets,
even if the topos does not have enough points, see \cite[Section ``For logical reasons'']{Bek04}.
\end{remark}

Since $(\tau_{\le n}X)_m$ and $(\tau_{< n}X)_m$ might fail to be representable sheaves,
we can't make sense of our Kan condition, and thus can't say much about them at this level of generality.
But in our situation of interest $\cC=\{$Banach manifolds$\}$,
we will nevertheless be able to analyze the structure of $\tau_{\le n}X$ and $\tau_{< n}X$.
Our arguments are not completely formal, so we defer them to section \ref{Integrating Lie n-algebras}.

\section{Banach manifolds}\label{sec:BM}

From now on, we work with the category of manifolds, by which we mean Banach manifolds,
and infinitely differentiable maps 
in the sense of Fr\'echet \cite[Section I.3]{Lan72}.
Given two manifolds $M$, $N$, a map $f:M\to N$ is called a {\em submersion} if the induced maps on tangent spaces is split surjective.
This condition implies that $f$ is locally diffeomorphic to a projection $V\oplus W\to V$ \cite[Corollary I.5.2]{Lan72}.
Our working definition will be the following:

\begin{definition}
A map $f:M\to N$ is a submersion if for every point $x\in M$, 
there exists a neighborhood $U$ of $f(x)$ and a section $U\to M$ sending $f(x)$ to $x$.
\end{definition}

\noindent In this section, we shall prove that the category of Banach manifold, 
equipped with the Grothendieck pretopology in which the covers are the surjective submersions, satisfies Assumptions \ref{assu}.

\begin{remark}\rm
For the constructions of sections \ref{SEC5}, \ref{SEC6}, \ref{Integrating Lie n-algebras},
it might seem more natural to use Fr\'echet instead of Banach manifolds.
But the implicit function theorem fails to be true for Fr\'echet manifolds, and we don't know whether 
surjective submersions form a Grothendieck pretopology in that context.
\end{remark}

We recall a couple basic facts about transversality in Banach manifolds (see \cite[Section I.5]{Lan72}, \cite[Section 3.5]{AMR88} for these and other related results).

\begin{lemma}\label{tfl}
Let $f:M\to N$ be a submersion and $y\in N$ be a point.
Then $f^{-1}(y)$ is a manifold.
\end{lemma}

\proof
Let $x\in f^{-1}(y)$ be a point, and let us identify $M$ locally with its tangent space $V:=T_x M$.
Let $s:N\to V$ be a local section of $f$.
We then have a direct sum decomposition $V=\ker(T_x f)\oplus\im(T_y s)$.
Let $p:V\to \ker(T_xf)$ denote the projection.

The derivative of $\phi:=(s\circ f) + p$ is the identity of $V$,
so by the inverse function theorem \cite[Theorem I.5.1]{Lan72} $\phi$ a local diffeomorphism.
The result follows since $f^{-1}(y)=\phi^{-1}(\ker (T_x f))$.
\qed

\begin{lemma}\label{lp}
Let $f:M\to N$ be a (surjective) submersion and let $g:Q\to N$ be an arbitrary map.
Then the pullback $P:=M\times_N Q$ is a manifold and its map to $Q$ is a (surjective) submersion.
\end{lemma}

\proof
The question being local, we may assume that $N$ is a vector space.
The space $P$ is then the preimage of the origin under the submersion $(f,-g):M\times Q\to N$ and is a manifold by Lemma \ref{tfl}.
The fact that $P\to Q$ is a (surjective) submersion is easy and left to the reader.
\qed

\begin{corollary}
Surjective submersions form a Grothendieck pretopology on the category of Banach manifolds.
\end{corollary}

We now show that the category of Banach manifolds satisfies the Assumptions \ref{assu}.
The first of those assumptions is trivial.
The fact that  the category of Banach manifolds is Cauchy complete is the content of the following lemma:

\begin{lemma}\label{etr}\dontshow{etr}
Let $M$ be a manifold and $p:M\to M$ be an idempotent map, then $p(M)$ is a manifold.
\end{lemma}

\proof
Let $x\in p(M)$ be a point, and let us identify $M$ locally with the tangent space $V:=T_xM$.
The linear map $T_xp:V\to V$ being idempotent, the vector space $V$ splits as
\[
V=\im(T_xp)\oplus\ker(T_xp).
\]
Let $q$ denote the linear projection $q:V\to\ker(T_xp)$, with kernel $\im(T_xp)$.
The map $\phi:=p+q$ is smooth and satisfies $T_x\phi=1_V$. 
So by the inverse function theorem, we conclude that $\phi$ is a diffeomorphism in a neighborhood of $x$.
The claim follows since locally $p(M)=\phi(\im(T_xp))$, and $\im(T_xp)$ is a vector space.
\qed

\begin{lemma}
The retract of a surjective submersion is a surjective submersion.
\end{lemma}

\proof
The retract of a surjective map is clearly surjective.
So we consider a commutative diagram like (\ref{retD})
\[
\xymatrix{
A\ar[r]_i\ar[d]^f&X\ar@/_8pt/[l]_r\ar[d]^g\\
B\ar[r]_i&Y\ar@/_8pt/[l]_r
}
\] 
and assume that $g$ is a submersion. Let $a\in A$ be a point.
Since $g$ is a submersion, we can find a neighborhood $U\subset Y$ of $y:=gi(a)$, and a local section $s:U\to X$ of $g$ sending $y$ to $i(a)$.
The preimage $V:=r^{-1}(U)$ is then a neighborhood of $b:=f(a)$ and $r\circ s\circ i:V\to A$ is a local section of $f$ sending $b$ to $a$.
This shows that $f$ is a submersion.
\qed

Finally, the last of Assumptions \ref{assu} follows from the fact that surjective submersions can be refined by usual open covers,
and that the `open covers' Grothendieck pretopology is subcanonical.

\begin{remark}\rm
An open cover being a special case of a surjective submersion, 
it follows that each one of the `surjective submersions' and the `open covers' Grothendieck pretopologies refines the other one.
A functor $\{\text{Banach manifolds}\}^{\mathit{op}}\!\to\{\text{Sets}\}$ is therefore 
a sheaf for one of the two pretopologies if and only if it is a sheaf for the other.
\end{remark}

We finish this section with an easy lemma that will be needed in Section \ref{Integrating Lie n-algebras}.

\begin{lemma}\label{aA4}
Let $f:M\to N$, $g:N\to Q$ be maps between Banach manifolds, with the properties that $f$ is surjective, and $g\circ f$ is a surjective submersion.
Then $g$ is a surjective submersion.
\end{lemma}

\proof
The surjectivity of $g$ is clear. To show that $g$ is a submersion, we pick a point $n\in N$
and construct a local section $Q\supset U\to N$, sending $g(n)$ to $n$.
Letting $m\in M$ be any preimage of $n$, we may pick a section $s:U\to M$ of $g\circ f$ sending $gf(m)=g(n)$ to $m$.
The map $f\circ s$ is the desired section of $g$.
\qed

\section{Integrating $L_\infty$-algebras}\label{SEC5}

In the coming sections, all $L_\infty$-algebras will be assumed to be non-negatively graded, and finite dimensional in each degree.

Given a nilpotent\footnote{
We use `nilpotent' as in \cite[Section 2]{Sul77}; it should not to be confused with the usual use of the term in commutative algebra.
} differential graded algebra $A$ over $\mathbb Q$, 
Sullivan \cite{Sul77} introduced its spatial realization
$\langle A\rangle$ (see also Bousfield-Guggenheim \cite{BG76}, in which $\langle A\rangle$ is called $FA$). 
It is a simplicial set with uniquely divisible homotopy groups, and nilpotent $\pi_1$.
The $n$-simplices of $\langle A\rangle$ are given by the set of $DGA$ homomorphisms from $A$ to the algebra $\Omega_{pol}^*(\Delta^n;\mathbb Q)$ 
of polynomial de Rham forms on the $n$-simplex with rational coefficients.
The idea of applying the above construction to the Chevalley-Eilenberg complex of an $L_\infty$-algebra is due to Getzler \cite{Get04}, 
who was inspired by Hinich \cite{Hin97}.

We shall be working with algebras over $\R$, and without any nilpotence assumption.
As suggested in \cite{Sul77}, we use an algebra which is bigger than $\Omega_{pol}^*$, and more suitable for this analytical context.
One of our results (Theorem \ref{tgo}) states the existence of natural manifold structures on the above sets of algebra homomorphisms.

As explained in Section \ref{sec:BM}, the failure of the implicit function theorem for Fr\'echet manifolds forces us to work with Banach manifolds instead.
This leads us to the following slightly inelegant convention.
Fix an integer $r\ge 1$.
Given a compact finite dimensional manifold $M$, possibly with boundary and corners, let $\Omega^*_{C^r}(M)$ be the differential graded algebra given by
\[
\Omega^n_{C^r}(M):=\{\alpha\in \Gamma_{C^r}(M;\Lambda^n T^*M)\,|\,d\alpha\in \Gamma_{C^r}(M;\Lambda^{n+1} T^*M)\},
\]
where $\Gamma_{C^r}(M,\xi)$ denotes the space of $r$ times continuously differentiable sections of a smooth vector bundle $\xi$.
The choice of $r$ is immaterial, and will often be suppressed from the notation.

\begin{convention}\label{Bc}\dontshow{Bc}
\rm Thereafter, $\Omega^n$ will always mean $\Omega_{C^r}^n$.
\end{convention}

The only properties of $\Omega_{C^r}^*(M)$ that we shall use are that it is a differential graded Banach algebra, 
and that it behaves well under restriction. Any other completion of $\Omega^*_{C^\infty}(M)$ with the above properties could be used instead.

\begin{definition}\label{dig}\dontshow{dig}
Let $L$ be an $L_\infty$-algebra, and let $C^*(L)$ be its Chevalley-Eilenberg complex (\ref{cocc}).
The integrating simplicial manifold $\intinf L$ is then given by \dontshow{in}
\begin{equation}\label{in}
\big(\intinf L\big)_m:=Hom_{DGA}\big(C^*(L),\Omega^*(\Delta^m)\big),
\end{equation}
where $\Omega^*(\Delta^m)$ denotes the differential graded algebra of de Rham forms on the $m$-simplex, following Convention \ref{Bc}.

Given a morphism $f:[m]\to[n]$ in $\Delta$, the corresponding map $(\intinf L)_n\to (\intinf L)_m$ is induced 
by the restriction $\Omega^*(\Delta^n)\to \Omega^*(\Delta^m)$ along $f_*\!:\!\Delta^m\to \Delta^n$.
\end{definition}

\begin{example}\rm\label{xls}\dontshow{xls}
Let $L=L_{n-1}$ be an $L_\infty$-algebra concentrated in degree $n-1$, and with all brackets vanishing.
Then $(\intinf L)_m=\Omega_{closed}^n(\Delta^m;L_{n-1})$.
\end{example}

\begin{example}\rm\label{xlw}\dontshow{xlw}
Let $L=L_n\oplus L_{n-1}$ be a contractible $L_\infty$-algebra concentrated in degrees $n$ and $n-1$. 
All its brackets are zero except $[\,]:L_n\to L_{n-1}$, which is an isomorphism.
Then $(\intinf L)_m=\Omega^n(\Delta^m;L_n)$.
\end{example}

\begin{example}\rm\label{xl}\dontshow{xl} {\cite[``Theorem''$(8.1)'$]{Sul77}}
Let $\g$ be a Lie algebra, viewed as $L_\infty$-algebra concentrated in degree $0$,
and let $G$ be a Lie group with Lie algebra $\g$.
A homomorphism of graded algebras $C^*(\g)\to\Omega^*(\Delta^m)$ is uniquely determined by a linear map
$\g^*\to \Omega^1(\Delta^m)$. 
This is the same as a $\g$-valued 1-form $\alpha\in\Omega^1(\Delta^m;\g)$, which
can then be interpreted as a connection on the trivial $G$-bundle $G\times \Delta^m\to\Delta^m$.
This graded algebra homomorphism respects the differentials if and only if $\alpha$ satisfies the Maurer-Cartan\footnote{
The usual Maurer-Cartan equation is $d\alpha+{\frac{1}{2}}[\alpha\wedge\alpha]=0$. One can bring (\ref{Mmc}) in that form by replacing $\alpha$ by $-\alpha$.
}
equation \dontshow{Mmc}
\begin{equation}\label{Mmc}
d\alpha={\frac{1}{2}}[\alpha\wedge\alpha],
\end{equation}
where $[\alpha\wedge\alpha]$ denotes the image under the bracket map of the $(\g\otimes\g)$-valued $2$-form $\alpha\wedge\alpha$.
Such $\alpha$'s correspond to flat connections on $G\times \Delta^m\to\Delta^m$.
Assigning to such a connection the set of all its horizontal sections establishes an isomorphism between the space of flat connections,
and the space of maps $\Delta^m \to G$, modulo translation.
Note that $\alpha$ being in $\Omega_{C^r}^1(\Delta^m;\g)$ is equivalent to that map being of class $C^{r+1}$.
So we get that \dontshow{yig}
\begin{equation}\label{yig}
\left(\intinf \mathfrak g\right)_m = Map\big(\Delta^m,G\big)/G.
\end{equation}
where $Map$ denote $C^{r+1}$ maps.
\end{example}

The main goal of this section is to prove Theorem \ref{tgo}, which states that $\intinf L$ is a Kan simplicial manifold.
For this purpose, we introduce a filtration (\ref{pwt}) whose associated graded pieces are the $L_\infty$-algebras of Examples 
\ref{xls}, \ref{xlw}, and \ref{xl}.

Given an $L_\infty$-algebra $L$, let $\partial_n:L_n\to L_{n-1}$ denote the components of the $1$-ary bracket $\partial=[\,]$.
The homology groups $H_n(L)=\ker(\partial_n)/\im(\partial_{n+1})$
behave very much like homotopy groups,
with $H_{n-1}$ playing the role of $\pi_n$.
This justifies our following terminology:

\begin{definition}\label{dhv}\dontshow{dhv}
The Postnikov pieces $\tau_{\le n}L$ and $\tau_{<n}L$ of an $L_\infty$-algebra $L$ are given by
\begin{equation*}
(\tau_{\le n}L)_i=\begin{cases}
L_i&\text{ if }i<n\\
\coker(\partial_{n+1})&\text{ if }i=n\\
0&\text{ if }i>n
\end{cases}
\qquad
(\tau_{< n}L)_i=\begin{cases}
L_i&\text{ if }i<n\\
\im(\partial_n)&\text{ if }i=n\\
0&\text{ if }i>n
\end{cases}
\end{equation*}
will all the brackets inherited from $L$.
\end{definition}

We then have a tower of projections \dontshow{pwt}
\begin{equation}\label{pwt}
\ldots \to \tau_{<n+1}L \to \tau_{\le n}L \to \tau_{<n}L \to \ldots
\to \tau_{\le 1}L\to 0
\end{equation}
whose inverse limit is $L$. The $i$th homology group of $\tau_{<n}L$ ($\tau_{\le n}L$) is that of $L$ for $i<n$ ($i\le n$) 
and zero otherwise, and the maps $ \tau_{<n}L \rightarrow\tau_{\le n-1}L$ are isomorphisms on homology.

Denote by $\Lambda^{m,j}\subset\Delta^m$ the geometric realization of $\Lambda[m,j]\subset\Delta[m]$.

\begin{lemma}\label{mg}\dontshow{mg}
Let $M$ be a smooth finite dimensional manifold, and $X$ be the simplicial manifold given by 
$X_m=Map(\Delta^m,M)$, where $Map$ denotes $C^{r+1}$ maps.
Then $X$ satisfies the Kan condition.
\end{lemma}

\proof
A map $\Lambda[m,j]\to X$ is the same thing as a map $\Lambda^{m,j}\to M$ which is $C^{r+1}$ when restricted to each face of $\Lambda^{m,j}$.

To show that (\ref{dav}) is surjective, we must show that every such map $f:\Lambda^{m,j}\to M$ extends to $\Delta^m$.
Pick a closed embedding $\iota:M\hookrightarrow \R^N$,
let $\nu(M)\subset\R^N$ be a tubular neighborhood of $M$, and let $\pi:\nu(M)\to M$ be the normal projection.
As before, let $[m]= \{0,1,\ldots, m\}$.
For $\{j\}\subsetneq I\subset[m]$, 
let $p_I:\Delta^m\to \Lambda^{m,j}$ denote the affine projection
given on vertices by \dontshow{defpI}
\begin{equation}\label{defpI}
p_I(e_i):=\begin{cases} e_j &\text{if}\, i\in I\\ 
e_i &\text{if} \, i\not\in I.
\end{cases}
\end{equation}
For $f:\Lambda^{m,j}\to M$, let $\tilde f:\Delta\to \R^N$ be the sum \dontshow{meqm}
\begin{equation}\label{meqm}
\tilde f:=\sum_{\{j\}\subsetneq I\subset[m]}
(-1)^{|I|}\,\iota\circ f\circ p_I.
\end{equation}
One easily checks that $\tilde f|_{\Lambda^{m,j}}=\iota\circ f$.
Given $\epsilon>0$, let $\Lambda_\epsilon^{m,j}\subset\Delta^m$ be the $\epsilon$\,-neighborhood of $\Lambda^{m,j}$.
and let $\epsilon(f):=\inf\{\epsilon:\tilde f|_{\Lambda^{m,j}_\epsilon}\subset\nu(M)\}$.
Pick smooth maps $r_\epsilon:\Delta^m\to\Lambda_\epsilon^{m,j}$ fixing $\Lambda^{m,j}$ and depending smoothly on $\epsilon$.
An extension $\bar f:\Delta^m\to M$ is then given by the formula
\[
\bar f:=\pi\circ\tilde f\circ r_{\epsilon(f)}.
\]
This provides a global lift of (\ref{dav}).

To show that (\ref{dav}) is submersive, we modify the above construction so that the image $\bar f_0$ of a chosen map $f_0$ is specified.
Instead of (\ref{meqm}), we set \dontshow{mek}
\begin{equation}\label{mek}
\tilde f':=
\sum(-1)^{|I|}\,\iota\circ f\circ p_I
+\big(\iota\circ \bar f_0-
\sum(-1)^{|I|}\,\iota\circ f_0\circ p_I\big).
\end{equation}
Let us now further assume that $r_\epsilon$ has been chosen so that $r_\epsilon=1$ if $\Lambda_\epsilon^{m,j}=\Delta^m$.
Letting $\epsilon'(f):=\inf\{\epsilon:\tilde f'|_{\Lambda^{m,j}_\epsilon}\subset\nu(M)\}$, we now have a new extension
\[
\bar f':=\pi\circ\tilde f'\circ r_{\epsilon'(f)},
\]
which is compatible with our initial choice $\bar f_0$.
\qed

\begin{lemma}\label{pfa}\dontshow{pfa}
The quotient of a Kan simplicial manifold $X$ by the proper free action of a Lie group $G$ is again a Kan simplicial manifold.
\end{lemma}

\proof
Since $\Lambda[m,j]$ is simply connected, we have
\[
Hom(\Lambda[m,j],X/G)=Hom(\Lambda[m,j],X)/G.
\]
The result then follows since the quotient of a surjective submersion by a proper free action is again a surjective submersion.
\qed

\begin{lemma}\label{ddt}\dontshow{ddt} 
The simplicial manifold $X$ given by $X_m=\Omega^n(\Delta^m)$ is Kan, and so is the one given by $X_m=\Omega^n_{closed}(\Delta^m)$.
\end{lemma}

\proof
Given a (closed) form $\alpha$ on $\Lambda^{m,j}$, we can extend it to $\Delta^m$ by the formula
\begin{equation}\label{swa}
\bar\alpha:=\sum_{\{j\}\subsetneq I\subset[m]}
(-1)^{|I|}\,p_I^*\alpha.
\end{equation}
where $p_I$ are the projections introduced in the proof of Lemma \ref{mg}. This shows that (\ref{dav}) is surjective.
If we have a chosen extension $\bar\alpha_0 \in\Omega^n(\Delta^m)$ of a form $\alpha_0\in\Omega^n(\Lambda^{m,j})$, 
we can modify (\ref{swa}) in a way similar to (\ref{mek})
\begin{equation*}
\bar\alpha':=
\sum(-1)^{|I|}\,p_I^*\alpha
+\big(
\bar\alpha_0
-\sum(-1)^{|I|}\,p_I^*\alpha_0
\big)
\end{equation*}
so that the extension of $\alpha_0$ is exactly $\bar\alpha_0$.
This shows that (\ref{dav}) is a submersion.
\qed

We are now in position to prove our main theorem.

\begin{theorem}\label{tgo}\dontshow{tgo}
Let $L$ be an $L_\infty$-algebra.
Then each $(\intinf L)_m$ is a smooth manifold, and $\intinf L$ satisfies the Kan condition (\ref{dav}).
Moreover, all the maps \dontshow{ptV}
\begin{equation}\label{ptV}
\ldots\,\rightarrow\, \intinf\tau_{\le 1}L\,\rightarrow\, \intinf\tau_{<1}L\,\rightarrow \, \intinf\tau_{\le 0}L\,\rightarrow\,*
\end{equation}
are Kan fibrations.
\end{theorem}

\proof
By Lemma \ref{cmg}, the maps (\ref{ptV}) being Kan fibrations implies that each $\intinf \tau_{<n} L$ is Kan.
Since $(\intinf L)_m=(\intinf\tau_{<n}L)_m$ for $n\ge m$, this also implies the Kan condition for $\intinf L$.

We first note that $\intinf\tau_{\le 0} L$ is Kan. 
Indeed, $\tau_{\le 0} L$ is an ordinary Lie algebra. 
Its integrating simplicial manifold has been computed in Example \ref{xl} and it is Kan by Lemmas \ref{mg} and \ref{pfa}.

To see that $\intinf\tau_{< n} L\to\intinf\tau_{\le n-1} L$ is a Kan fibration,
note that $C^*(\tau_{< n} L)$ is freely generated as a differential graded algebra over $C^*(\tau_{\le n-1} L)$
by the vector space $im(\partial_n)^*$, put in degree $n$.
It follows that  $\intinf\tau_{< n} L$ is the product of $\intinf\tau_{\le n-1} L$ with the simplicial manifold of Example \ref{xlw}.
Assuming by induction that the sheaves $(\intinf\tau_{< n} L)_m$ are manifolds, it follows that so are the sheaves $(\intinf\tau_{\le n-1} L)_m$.
By Lemmas \ref{rmf} and \ref{ddt}, we also deduce that $\intinf\tau_{< n} L\to\intinf\tau_{\le n-1} L$ is a Kan fibration.

We now show that \dontshow{kqz}
\begin{equation}\label{kqz}
\intinf\tau_{\le n-1} L\longrightarrow\intinf\tau_{< n-1} L
\end{equation}
is a Kan fibration.
Let us use the shorthand notations $H:=H_{n-1}(L)$ and $\g:=H_0(L)$. 
Pick a splitting of the inclusion $H\to (\tau_{\le n-1}L)_{n-1}$.
The graded algebra $C^*(\tau_{\le n-1}L)$ is freely generated over $C^*(\tau_{< n-1}L)$ by the 
vector space $H^*$.
A graded algebra homomorphism $\varphi:C^*(\tau_{\le n-1}L)\to\Omega^*(\Delta^k)$ is therefore the same thing as 
a homomorphism  $\psi:C^*(\tau_{< n-1}L)\to\Omega^*(\Delta^k)$ along with a $H$-valued $n$-form $\beta$.
We now characterize the pairs $(\psi,\beta)$ such that $\varphi=\varphi(\psi,\beta)$ respects the differential.
Clearly, $\psi$ needs to respect the differential,
so we assume that this is the case and investigate the condition on $\beta$.
By the Leibnitz rule, it is enough to verify on $H^*$ whether $d\circ\varphi=\varphi\circ\delta$.

Let us write $\delta|_{H^*}$ as $\delta_0+\delta_1$, where 
$\delta_0$ lands in $L_0^*\otimes H^*$ and $\delta_1$ lands in $C^*(\tau_{< n-1}L)$.
Since $\partial=[\,]$ is a derivation of $[\,,]$ the component of $[\,,]$ that maps $im(\partial_1)\otimes H$ to $H$ vanishes.
It follows that $\delta_0$ only lands in $\g^*\otimes H^*$.
Let $\alpha\in\Omega^1(\Delta^k;\g)$ be the restriction of $\varphi$ to $\g^*$,
and $\gamma\in\Omega^{n+1}(\Delta^k;H)$ be the composite $\varphi\circ\delta_1$.
We then have $d\circ\varphi|_{H^*}=d\beta$ and 
$\varphi\circ\delta|_{H^*}=\varphi\circ\delta_0+\varphi\circ\delta_1=[\alpha\wedge\beta]+\gamma$.
It follows that $\varphi(\psi,\beta)$ respects the differentials iff $\beta$ satisfies \dontshow{abg}
\begin{equation}\label{abg}
d\beta=[\alpha\wedge\beta]+\gamma.
\end{equation}

Let $G$ be the simply connected Lie group integrating $\g$.
Since $H$ is a $\g$-module, we can integrate it to a $G$-module.
The 1-form $\alpha$ satisfies the Maurer-Cartan equation $d\alpha=\frac{1}{2}[\alpha,\alpha]$ so
we can integrate it to a map $f:\Delta^k\to G$ satisfying $-\alpha=f^{-1}df$.
Let's also fix the image of some vertex, so as to make $f$ uniquely defined.
Let now \dontshow{bg}
\begin{equation}\label{bg}
\beta':=f\beta \quad\qquad\gamma':=f\gamma.
\end{equation}
Clearly, once $\psi$ is fixed, the data of $\beta$ and $\beta'$ are equivalent to each other.
But
in terms of $\beta'$, equation (\ref{abg}) becomes much simpler: \dontshow{jwh}
\begin{equation}\label{jwh}
d\beta'=\gamma'.
\end{equation}
We want to solve (\ref{jwh}) for $\beta'$, so we first check that $\gamma'$ is closed. Indeed, we have
\[
d\gamma'=[dff^{-1}\wedge f\gamma]+f d\gamma
=-[f\alpha\wedge f\gamma]+f d\gamma
=f\big(d\gamma-[\alpha\wedge\gamma]\big)=0
\]
since
\[
\begin{split}
d\gamma&=d(d\beta-[\alpha\wedge\beta])=[\alpha\wedge d\beta]-[d\alpha\wedge\beta]\\
&={\textstyle[\alpha\wedge\gamma]+[\alpha\wedge[\alpha\wedge\beta]]-{\frac{1}{2}}}[[\alpha\wedge\alpha]\wedge\beta]=[\alpha\wedge\gamma],
\end{split}
\]
where the last equality holds because $[\,,]:\g\otimes H\to H$ is a Lie algebra action.

Going back to our problem of showing that (\ref{kqz}) is a Kan fibration,
we must show that \dontshow{jw}
\begin{equation}\label{jw}
(\intinf\tau_{\le n-1} L)_m\,\longrightarrow\, Hom\big(\Lambda[m,j]\to\Delta[m],\,\intinf\tau_{\le n-1} L\to \intinf\tau_{< n-1} L\big)
\end{equation}
is a surjective submersion.
An element in the RHS of (\ref{jw}) consists of two compatible homomorphisms $\psi:C^*(\tau_{< n-1} L)\to\Omega^*(\Delta^m)$ and
$\varphi:C^*(\tau_{\le n-1} L)\to\Omega^*(\Lambda^{m,j})$.
Let $P$ be the pullback of
\begin{equation*}
\Omega^n(\Lambda^{m,j})\stackrel{d}{\longrightarrow}\Omega^{n+1}_{\text{\it closed}}(\Lambda^{m,j})
\longleftarrow\Omega^{n+1}_{\text{\it closed}}(\Delta^m),
\end{equation*}
and let $\Omega^n(\Delta^m)\to P$ be the obvious map.
We then have a commutative square \dontshow{asp}
\begin{equation}\label{asp}
\begin{matrix}
\xymatrix{
(\intinf\tau_{\le n-1} L)_m\ar[r]\ar[d]&\Omega^n(\Delta^m)\ar[d]\\
Hom\big(\Lambda[m,j]\to\Delta[m],\intinf\tau_{\le n-1} L\to \intinf\tau_{< n-1} L\big)\ar[r]&P.}
\end{matrix}
\end{equation}
The bottom map assigns to a pair $(\varphi,\psi)$ the forms $\beta'\in\Omega^n(\Lambda^{m,j})$, 
$\gamma'\in\Omega^{n+1}_{closed}(\Delta^m)$ defined in (\ref{bg}).
And the top map assigns to $\bar\varphi:C^*(\tau_{\le n-1} L)\to\Omega^*(\Delta^m)$ the corresponding form $\bar\beta'\in\Omega^n(\Delta^m)$.
By the previous discussion, an extension $\bar\varphi$ of $\varphi$ compatible with a given $\psi$, is equivalent to an 
extension $\bar\beta'$ of $\beta'$ compatible with a given $\gamma'$.
In other words, (\ref{asp}) is a pullback square.
The right vertical map is a surjective submersion by Lemma \ref{dpo}, therefore so is the left vertical map.

Lastly, we show by induction on $m$ that $(\intinf\tau_{\le n-1} L)_m$ is a manifold.
The space on the lower left corner of (\ref{asp}) is a manifold by Corollary \ref{ikq}. 
The result then follows from Lemma \ref{dpo}.
\qed

\begin{lemma}\label{dpo}\dontshow{dpo}
Let $P$ be the pullback of
\begin{equation*}
\Omega^n(\Lambda^{m,j})\stackrel{d}{\longrightarrow}\Omega^{n+1}_{closed}(\Lambda^{m,j})\longleftarrow\Omega^{n+1}_{closed}(\Delta^m).
\end{equation*}
Then the natural map $\Omega^n(\Delta^m)\to P$ is a surjective submersion.
\end{lemma}

\proof
We first show surjectivity.
Let $\alpha\in \Omega^n(\Lambda^{m,j})$ , $\beta\in \Omega^{n+1}_{closed}(\Delta^m)$ be forms
satisfying $d\alpha=\beta|_{\Lambda^{m,j}}$.
We must find an extension $\bar\alpha\in\Omega^n(\Delta^m)$ of $\alpha$ satisfying $d\bar\alpha=\beta$.

Let $p_I:\Delta^m\to\Lambda^{m,j}$ be as in (\ref{defpI}). 
For $s\in[0,1]$, let $p_I^s:=s\cdot Id+(1-s)p_I$, and let $v_I$ be the vector field given by $(v_I)_x=x-p_I(x)$.
The form \dontshow{stj}
\begin{equation}\label{stj}
\bar\alpha_I:=p_I^*\alpha+\int_{s=0}^1(p_I^s)^*\iota_{v_I}\beta \,ds
\end{equation}
then satisfies $\bar\alpha_I\big|_{p_I(\Delta^m)}=\alpha$, and $$\bar\alpha_J\big|_{p_I(\Delta^m)}=\bar\alpha_{J'}\big|_{p_I(\Delta^m)}$$ when 
$p_I(\Delta^m)\cap p_J(\Delta^m)=p_I(\Delta^m)\cap p_{J'}(\Delta^m)$.
Its differential is given by \dontshow{baz}
\begin{equation}\label{baz}
d\bar\alpha_I=p_I^*d\alpha+\int_{s=0}^1(p_I^s)^*\mathcal L_{v_I}\beta \,ds.
\end{equation}
where $\mathcal L$ denotes the Lie derivative.
Since $\beta$ is closed and $p_I^s$ is the flow generated by $v_I$, this becomes \dontshow{mes}
\begin{equation}\label{mes}
d\bar\alpha_I=p_I^*d\alpha+(p_I^1)^*\beta-(p_I^0)^*\beta=\beta.
\end{equation}
Now let $\bar\alpha$ be given by \dontshow{wuv}
\begin{equation}\label{wuv}
\bar\alpha=\sum_{\{j\}\subsetneq I\subset[m]}(-1)^{|I|}\bar\alpha_I.
\end{equation}
By inclusion exclusion, we have $\bar\alpha|_{\Lambda^{m,j}}=\alpha$, and by (\ref{mes}) it satisfies $d\bar\alpha=\beta$.

We now show that $\Omega^n(\Delta^m)\to P$ is a submersion.
Let $(\alpha^0,\beta^0)\in P$, and let $\bar\alpha^0\in\Omega^n(\Delta^m)$ be a preimage.
We modify (\ref{wuv}) to produce a section $P\to \Omega^n(\Delta^m)$ sending $(\alpha^0,\beta^0)$ to $\bar\alpha^0$.
This can again be done by an explicit formula:
\begin{equation*}
(\alpha,\beta)\,\mapsto\,\bar\alpha^0+
\sum(-1)^{|I|}
(\bar\alpha_I-\bar\alpha_I^0).
\end{equation*}\qed

\section{The homotopy groups of $\intinf L$}\label{SEC6}

From a purely homotopy theoretic point of view, the simplicial manifold $\intinf L$ is not very interesting.
Each manifold $(\intinf L)_m$ is contractible, and so is the geometric realization of $\intinf L$.
Indeed, a nullhomotopy $h_t:\Delta^m\to\Delta^m$, $t\in[0,1]$, induces nullhomotopies
$h_t^*:\Omega^*(\Delta^m)\to\Omega^*(\Delta^m)$, and
\[
h_t^*\circ-:Hom_{DGA}\big(C^*(L),\Omega^*(\Delta^m)\big)\to Hom_{DGA}\big(C^*(L),\Omega^*(\Delta^m)\big),
\]
proving that each $(\intinf L)_m$ is contractible.

The simplicial homotopy groups of $\intinf L$ carry much more information.
In general, these simplicial homotopy groups can fail to be representable sheaves (see Example \ref{eND}).
But the worse that can happen is that $\pi_n^\simpl(\intinf L)$ is the quotient of a Lie group by a non-discrete finitely generated subgroup (Theorem \ref{tls}).
In that sense, it is never very far from being a Lie group.

If $X$ is a simplicial manifold, then the space
\begin{equation}\label{PI}
\hspace{3cm}\big\{x\in X_n\,\big|\, d_i(x)=*\, \forall i\big\}\big/\sim\qquad (\,\text{\small $\sim$ as in (\ref{1dpi})}\hspace{.02cm})
\end{equation}
is the best possible approximation of $\pi_n^{\simpl}(X)$ by a topological space. 
If we assume that $\{x\in X_n| d_i(x)=*\, \forall i\}$ is a submanifold of $X_n$,
that $\sim$ is a Hausdorff equivalence relation given by a foliation,
and that the space
\[
\big\{x_1,x_2\in X_n\,\big|\, d_i(x_1)=d_i(x_2)=*,\, x_1\sim x_2\big\}
\]
admits local lifts into $X_{n+1}$,
then (\ref{PI}) is actually a manifold, and that manifold represents $\pi_n^{\simpl}(X)$.
 
Below are some examples of simplicial manifolds whose simplicial homotopy groups are representable.
We just check that the spaces (\ref{PI}) are Lie groups.

\begin{example}\rm\label{xig}\dontshow{xig}
Let $G$ be a Lie group with Lie algebra $\g$. 
Then $\pi_1^\simpl(\intinf\g)$ is the universal cover of $G$, while $\pi_n^\simpl(\intinf\g)=\pi_n(G)$ for $n\ge 2$.

Indeed, an element of $\pi_n^\simpl(\intinf\g)$ is represented by a map 
$f:\Delta^n\to G$, which is constant on each facet of $\Delta^n$.
Two maps are identified if one is a translate of the other, or if they are homotopic relatively to $\partial\Delta^n$. 
For $n\ge 2$, the above condition means that $f|_{\partial\Delta^n}$ is constant.
There's a unique translate of $f$ sending $\partial\Delta^n$ to $e\in G$, and we recover the usual definition of $\pi_n(G)$.
For $n=1$, there's no condition on $f$. 
Again, there's a unique translate sending a given vertex to $e$.
The group $\pi_1^\simpl(\intinf\g)$ is then the based paths space of $G$ modulo homotopy fixing the endpoints, i.e.
the universal cover of $G$.
\end{example}

\begin{example}\rm\label{xah}\dontshow{xah}
Let $L$ be as in Example \ref{xls}. 
Then $\pi_n^\simpl(\intinf L)=L_{n-1}$, and all other homotopy groups vanish.

An element of $\pi_m^\simpl(\intinf L)$ is represented by a closed $L_{n-1}\hspace{.02cm}$-valued $n$-form $\omega$ on 
$\Delta^m$ that vanishes when restricted to $\partial\Delta^m$.
Two forms $\omega$, $\omega'$ are then equivalent if
$\omega\sqcup\omega'$, defined on the $m$-sphere $S^m=\Delta^m\sqcup_{\partial\Delta^m}\Delta^m$, 
extends to a closed form $\bar\omega$ on the $(m+1)$-disk.

If $m\not = n$, then any two forms are equivalent. 
Indeed, since $\omega\sqcup\omega'$ is zero in $H^n(S^m;L_{n-1})$, we may pick a form $\alpha$ such that $d\alpha=\omega$.
Using formulas like (\ref{swa}) and partitions of unity, we can extend $\alpha$ to a form $\bar\alpha$ on the $(m+1)$-disk.
The desired extension is then given by $\bar\omega=d\bar\alpha$.
This shows that $\pi_m(\intinf L)=0$ for $m\not = n$.

If $m=n$, then by Stokes' theorem any extension $\bar\omega$ must satisfy \dontshow{stl}
\begin{equation}\label{stl}
\int_{D^{n+1}}d\bar\omega=\int_{\Delta^n}\omega-\int_{\Delta^n}\omega'.
\end{equation}
A necessary condition for the existence of a closed extension $\bar\omega$ is for the right hand side of (\ref{stl}) to vanish.
If it does vanish, then the same argument as above shows that $\omega$ and $\omega'$ are equivalent.
We conclude that $\omega\sim\omega'$ iff theirs integrals agree, and hence that $\pi_n(\intinf L)=L_{n-1}$.
\end{example}

\begin{example}\rm\label{yaa}\dontshow{yaa}
Let $L$ be as in Example \ref{xlw}. 
Then all the simplicial homotopy groups of $\intinf L$ vanish.

As above, an element of $\pi_m^\simpl(\intinf L)$ is given by a form $\omega$ on 
$\Delta^m$, and two forms $\omega$, $\omega'$ are equivalent if
$\omega\sqcup\omega'$ extends to the $(m+1)$-disk.
There are no obstructions to finding such an extension, and thus any two forms are equivalent.
\end{example}

Even though $\pi_n^\simpl(\intinf L)$ might fail to be a Lie group, it still belongs to a class 
of group sheaves for which it makes sense to take the Lie algebra.
These are quotient sheaves of the form $G/A$, where $G$ is a finite dimensional Lie group, and $A$
is a finitely generated subgroup of $G$, possibly non-discrete.
We shall call them {\em finite dimensional diffeological groups}.
The universal cover of a finite dimensional diffeological group being a Lie group \cite{Sou80}, 
we can define its Lie algebra to be that of its universal cover.

\begin{theorem}\label{tls}\dontshow{tls}
Let $L$ be an $L_\infty$-algebra, and let $G$  be the simply connected Lie group integrating $H_0(L)$.
Then the simplicial homotopy groups $\pi_n^\simpl(\intinf L)$ are finite dimensional diffeological groups.
We have $\pi_1^\simpl(\intinf L)= G$,
and a long exact sequence \dontshow{lle}
\begin{equation}\label{lle}
\ldots\to\pi_{n+1}^\simpl(\intinf L)\to\pi_{n+1}(G)\to H_{n-1}(L)\to
\pi_n^\simpl(\intinf L)\to\pi_n(G)
\to\ldots
\end{equation}
for $n\ge 2$.
\end{theorem}

\proof
We work by induction on the tower (\ref{ptV}).
The initial case $L=\tau_{\le 0}L$ is the computation done in Example \ref{xig}.

Assuming that the theorem holds for $\tau_{\le n-1}L$, we show that it also holds for $\tau_{<n}L$. 
Indeed, the fiber of $\intinf\tau_{<n}L\to\intinf\tau_{\le n-1}L$ is the simplicial manifold of Example \ref{xlw}.
Its homotopy groups being zero by Example \ref{yaa}, we conclude by Proposition \ref{sess} that 
$\pi_m^\simpl(\intinf\tau_{< n}L)=\pi_m^\simpl(\intinf\tau_{\le n-1}L)$. 
Similarly, we have $\pi_m(\tau_{< n}L)=\pi_m(\tau_{\le n-1}L)$.
The sequence (\ref{lle}) is therefore the same as that for $\tau_{\le n-1}L$.

Now, we assume that the result holds for $\tau_{< n-1}L$ and
show that it also holds for $\tau_{\le n-1}L$. 
The fiber of $\intinf\tau_{\le n-1}L\to\intinf\tau_{< n-1}L$ is as in Example \ref{xls}, but with $H_{n-1}(L)$ instead of $L_{n-1}$. 
By Example \ref{xah}, its unique homotopy group is then $H_{n-1}(L)$ in dimension $n$.
It follows from Proposition \ref{sess} that $\pi_m^\simpl(\intinf\tau_{\le n-1}L)=\pi_m^\simpl(\intinf\tau_{< n-1}L)$ for $m\not =n,n+1$, 
and that we have an exact sequence
\[
\begin{split}
0\to\pi_{n+1}^\simpl(\intinf\tau_{\le n-1}L)\to\pi_{n+1}^\simpl(\intinf\tau_{< n-1}L)&\\
\to H_{n-1}(L)\to\,&\,\pi_{n}^\simpl(\intinf\tau_{\le n-1}L)\to\pi_{n}^\simpl(\intinf\tau_{< n-1}L)\to 0.
\end{split}
\]
By assumption, we also have $\pi_{n+1}^\simpl(\intinf\tau_{< n-1}L)=\pi_{n+1}(G)$.
A diagram chase shows the exactness of (\ref{lle}).

To see that $\pi_n^\simpl(\intinf L)$ is a finite dimensional diffeological group, consider the exact sequence
\[
0\,\to\, C\,\to\, H_{n-1}(L)\,\to\,\pi_n^\simpl(\intinf L)\,\to\, K \,\to\, 0\,,
\]
where $C$ and $K$ are the finitely generated groups $\coker(\pi_{n+1}^\simpl(\intinf L)\to \pi_{n+1}(G))$ and $\ker(\pi_n(G)\to H_{n-2}(L))$ respectively.
The sheaf $\pi_n^\simpl(\intinf L)$ is an extension of the discrete group $K$ by the finite dimensional diffeological group $H_{n-1}(L)/C$, and is thus itself a finite dimensional diffeological group.
\qed

An alternative, more conceptual proof of Theorem \ref{tls} can be given using the language of spectral sequences.
The Postnikov tower of $L$ induces a filtration of $\intinf L$.
The associated spectral sequence is then
\[
E^1_{m,n}=\pi^\simpl_m\big(\intinf \tau_{=n-1}L\big)\Rightarrow\pi_{m+n}^\simpl\big(\intinf L\big),
\]
where $\tau_{=n-1}L$ denotes the $L_\infty$-algebra with only $H_{n-1}(L)$ in degree $n-1$.
By the computations done in Examples \ref{xig} and \ref{xah}, the spectral sequence looks like this:

\[
\put(0,-83){\line(1,0){290}}
\put(0,-83){\line(0,1){95}}
\xymatrix@R.3cm@C.5cm{
0 & 0 & 0 & H_3(L) & 0 \\
0 & 0 & H_2(L) & 0 & 0 \\
0 & H_1(L) & 0 & 0 & 0 \\
G & \pi_2(G) & \pi_3(G)\ar[ul] & \pi_4(G\ar[uul]) & \pi_5(G)\ar[uuul] \\
\pi_1^\simpl(\intinf L) & \pi_2^\simpl(\intinf L) & \pi_3^\simpl(\intinf L) & \pi_4^\simpl(\intinf L) & \pi_5^\simpl(\intinf L) 
}
\]
Being very sparce, it can be reinterpreted as a long exact sequence.

\begin{corollary}
The Lie algebra of $\pi_n^\simpl(\intinf L)$ is canonically isomorphic to $H_{n-1}(L)$.
\end{corollary}

\proof
Apply the Lie algebra functor to the sequence (\ref{lle}).
\qed

\section{Integrating Lie $n$-algebras}\label{Integrating Lie n-algebras}

A Lie $n$-algebra is an $L_\infty$-algebra concentrated in degrees $\{0,\ldots,n-1\}$.
As explained before, the geometric realization of $\intinf L$ is contractible.
However, if $L$ is an Lie $n$-algebra, there is a modification of $\intinf L$ which has a more interesting homotopy type.
It is obtained by applying to $\intinf L$ the truncation functor $\tau_{\le n}$ of Definition \ref{DEFTAU}.

Following our previous convention, our Lie $n$-algebras will all be finite dimensional.

\begin{definition}\label{dism}\dontshow{dism}
Let $L$ be a Lie $n$-algebra.
Its integrating Lie $n$-group is then given by $\tau_{\le n}\intinf L$.
\end{definition}

It can happen that $(\tau_{\le n}\intinf L)_m$ fails to be a manifold (Example \ref{eND}),
in which case calling $\tau_{\le n}\intinf L$ a Lie $n$-group is an abuse of terminology.
But we know from Lemma \ref{check n-group} that the condition in Definition \ref{fk} is at least satisfied.

\begin{example}\rm
Let $\g$ be a Lie algebra, viewed as Lie 1-algebra, and let $G$ be the simply connected Lie group with Lie algebra $\g$.
Then its integrating Lie 1-group $\tau_{\le 1}\intinf \g$ is the simplicial manifold $K(G,1)$ described in (\ref{kg1}).

We just check it at the level of points; it is not much harder to check it at the level of sheaves.
We know from (\ref{yig}) that $(\intinf \g)_m=Map(\Delta^m,G)/G$.
Two points of $(\intinf \g)_m$ are identified in $(\tau_{\le 1}\intinf \g)_m$
if they are simplicially homotopic relatively to $\sk_0\Delta^m$.
The group $G$ being simply connected, this happens if any only if their restrictions to $\sk_0\Delta^m$ agree.
We thus get
\[
(\tau_{\le 1}\intinf \g)_m=Map\big(\sk_0\Delta^m,G\big)/G,
\]
which is just another notation for the $m$-th manifold in (\ref{kg1}).
\end{example}

\begin{remark}\rm
For $n=1$, we may generalize the above definition to Lie algebroids.
The Lie groupoid corresponding to $\tau_{\le 1}\intinf L$ is then exactly
the same as the one constructed by Crainic and Fernandez in \cite{CF03}.
In that same context, a variant of $\tau_{\le 2}\intinf L$ has been studied by Tseng and Zhu \cite{TZ06} 
using the formalism of (weak) groupoid objects in the category of smooth \'etale stacks.
\end{remark}

To better work with simplicial sheaves, we extend 
the Grothendieck pretopology of section \ref{sec:BM} to the category of sheaves on manifolds.

\begin{definition}\label{63}
Let $X$, $Y$ be sheaves on the category of Banach manifolds.
A map $X\to Y$ is called a submersion if for every Banach manifold $U$, point $u\in U$, element $y\in Y(U)$, and lift $q \in X(\{u\})$ of $y|_{\{u\}}$, 
there exists a neighborhood $U'\subset U$ of $u$, and a lift $x\in X(U')$ of $y|_{U'}$ such that $x|_{\{u\}}=q$.
\[
\xymatrix{
\{u\}\ar[d]\ar[r]^q&X\ar[d]\ar@{<.}[dl]+<-16pt,7pt>^(.41)x|!{[l];[dl]}{\phantom{\big|*\big|}}\\
U'\subset U\hspace{.8cm}&Y\ar@{<-}[l]+<7pt,0pt>_(.53)y
}
\]
A map is called a surjective submersion if it is surjective, and it is a submersion.
\end{definition}

Our next goal is the following result.

\begin{theorem}\label{TVF}\dontshow{TVF}
If $L$ is a Lie $n$-algebra, then the sheaf $(\tau_{\le n}\intinf L)_m$ is a manifold for every $m<n$.

Let $\partial_n:\pi_{n+1}(G)\to H_{n-1}(L)$ denote the boundary map appearing in the long exact sequence (\ref{lle}).
Then $(\tau_{\le n}\intinf L)_n$ is a manifold if and only if $\partial_n$ has discrete image.

Moreover, if $\partial_n$ has discrete image, 
then $(\tau_{\le n}\intinf L)_m$ is a manifold for every $m$, and the simplicial manifold $\tau_{\le n}\intinf L$ satisfies the Kan condition.
\end{theorem}

Given a simplicial manifold $X$ let us denote by $Z_n$ the sheaf $Hom(\partial\Delta[n+1],X)$, and by $B_n$ the image sheaf
$B_n:=\im\big(X_{n+1}\to Z_n\big)$.
Letting $S=\partial\Delta[n+1]$ in Definition \ref{S and SI}, we get an obvious surjection $p:Z_n\to\pi_n^{\simpl}(X)$.

\begin{lemma}\label{kz}\dontshow{kz}
Let $X$ be a reduced Kan simplicial manifold whose simplicial homotopy groups $\pi_m^{\simpl}(X)$ are finite dimensional diffeological groups for all $m\le n$.
Then $Z_n$ and $B_n$ are manifolds, and the map $p:Z_n\to\pi_n^\simpl(X)$ is a surjective submersion. 
\end{lemma}

\proof
We prove the lemma by induction on $n$.
For $n=0$, the result is trivial.
So let $n\ge 1$, and let us assume that $B_{n-1}$ is known to be a manifold.
Let $D[n]$ be the disk
\dontshow{Dis}
\begin{equation}\label{Dis}
D[n]:=\Delta[n]\cup(\partial\Delta[n]\times\Delta[1]),
\end{equation}
and note that $\partial D[n]=\partial\Delta[n]$.
By Lemma \ref{lmD}, the bottom arrow in the pullback
\[
\xymatrix{
Hom\big(D[n]/\partial D[n],X\big)\ar[d]\ar[r]\ar@{}[rd]|(.45){\line(0,1){5}\line(-1,0){5}}&\{*\}\ar[d]\\
Hom\big(D[n],X\big)\ar[r]&B_{n-1}
}
\]
is a submersion.
The sheaf $Hom(D[n],X)$ being a manifold by Corollary \ref{cor-lith},
it follows that $Hom(D[n]/\partial D[n],X)$ is a manifold.

Let $f:D[n]/\partial D[n]\to\partial\Delta[n+1]$ be the obvious weak equivalence, and let $C_f$ be its mapping cylinder.
Applying Corollary \ref{cor-lith} to the inclusion $D[n]/\partial D[n]\hookrightarrow C_f$,
we deduce that $Hom(C_f,X)$ is a manifold.
Since $\partial\Delta[n+1]$ is a retract of $C_f$, $Hom(\partial\Delta[n+1],X)$ is a retract of $Hom(C_f,X)$,
and so by Lemma \ref{etr}, it follows that $Z_n=Hom(\partial\Delta[n+1],X)$ is a manifold.

Let us use $\pi_n$ as shorthand notation for $\pi_n^{\simpl}(X)$.
We now show that $p:Z_n\to \pi_n$ satisfies the condition in Definition \ref{63}.
Let $U$ be a manifold, $u\in U$ be a point, and let us consider $\gamma\in\pi_n(U)$, $z\in Z_n$ satisfying $p(z)=\gamma|_{\{u\}}$.
Since $p$ is surjective, $\gamma$ lifts to an element $\tilde\gamma\in Z_n(U')$ for some neighborhood $U'\subset U$ of $u$.
The two maps
\[
z,\tilde\gamma|_{\{u\}}:\partial\Delta[n+1]\to X
\]
represent the same element in $\pi_n(\{u\})$, so there exists a simplicial homotopy 
$h:\partial\Delta[n+1]\times\Delta[1]\to X$ between them.
Let $H$ be the manifold $Hom(\partial\Delta[n+1]\times\Delta[1],X)$ and let $r_0, r_1:H\to Z_n$ be the two restrictions,
which are submersions by Corollary \ref{cor-lith}.
The element $h\in H$ maps to $\tilde\gamma|_{\{u\}}$ under $r_0$ and thus extends to a map $\theta:U''\to H$
satisfying $\theta(u)=h$, $r_0 \circ \theta=\tilde\gamma|_{U''}$, for some neighborhood $U''\subset U'$ of $u$.
The element
\[
\delta:= r_1\circ \theta\in Map(U'',Z_n)=Z_n(U'')
\]
then satisfies $\delta|_{\{u\}}=z$ and $p(\delta)=\gamma|_{U''}$. 

We now show that $B_n$ is a manifold.
By assumption, $\pi_n^\simpl(X)$ is of the form $\put(0,0){$\phantom{\pi}^\simpl$}\widetilde{\pi_n\hspace{.17cm}(X)}/A$, 
with 
$\put(0,0){$\phantom{\pi}^\simpl$}\widetilde{\pi_n\hspace{.17cm}(X)}$
a Lie group, and $A$ a finitely generated subgroup.
Consider the following pullback squares
\[
\xymatrix{
B_n\ar[d]\ar[r]\ar@{}[rd]|(.45){\line(0,1){5}\line(-1,0){5}}&\widetilde{Z_n}\ar[d]^{\widetilde p}\ar[r]\ar@{}[rd]|(.45){\line(0,1){5}\line(-1,0){5}}&Z_n\ar[d]^p\\
\{0\}\ar[r]&
\put(0,0){$\phantom{\pi}^\simpl$}\widetilde{\pi_n\hspace{.17cm}(X)}
\ar[r]&\pi_n^\simpl(X).
}
\]
The map $\widetilde{Z_n}\to Z_n$ is an $A$-principal bundle whose base is a manifold.
It follows that $\widetilde{Z_n}$ is a manifold.
The map $p$ is a submersion, and therefore so is $\widetilde p$.
By transversality, it follows that $B_n$ is a manifold.
\qed

\begin{lemma}\label{lmD}\dontshow{lmD}
Let $D[n]$ be as in (\ref{Dis}).
Then the restriction map
\[
r:Hom(D[n],X) \to B_{n-1}
\]
is a submersion.
\end{lemma}

\proof
We will prove that for each map $\gamma:U\to B_{n-1}$ and pair $x\in Hom(D[n],X)$, $u\in U$ such that $r(x)=\gamma(u)$,
there exists a local lift of $\gamma$ around $u$, that is, an open neighborhood $V\subset U$ of $u$, and a map
$\delta:V\to Hom(D[n],X)$ such that $r\circ \delta=\gamma$.

Since \(r':Hom(\Delta[n],X) \to B_{n-1}\) is a surjective map (of sheaves), 
there exists a neighborhood $U'\subset U$ of $u$, and a map $\gamma':U'\to Hom(\Delta[n],X)$ satisfying $r'\circ\gamma'=\gamma$.
The maps $x$ and $\gamma'(u)$ together define an element $a\in\pi_n^\simpl(X)$.
Let $\tilde\gamma:U''\to Hom(\Delta[n],X)$, $u\in U''\subset U'$, be a map representing the composite of $\gamma'$ and $a$.
The two maps $x:D[n]\to X$ and $\tilde\gamma(u):\Delta[n]\to X$ are now simplicially homotopic relatively to their boundary.

Let $H:=Hom(\Delta[n]\times\Delta[1],X)$ and $r_0:H\to Hom(\Delta[n],X)$, $r_1:H\to Hom(D[n],X)$ be the two restrictions maps.
We then get an element $h\in H$ such that $r_0(h)=\tilde\gamma(u)$ and $r_1(h)=x$. 
By Corollary \ref{cor-lith}, $r_0$ is a submersion,
so we can 
find a neighborhood $V\subset U''$ of $u$ and an extension
$\theta:V\to H$ of $h$ such that $r_0\circ\theta=\tilde\gamma$.
The map $r_1\circ \theta:V\to Hom(D[n],X)$ is the desired lift of $\gamma$.
\qed

The following is an easy variant of Corollary \ref{cor-lith}.

\begin{lemma}\label{KX}
Let $X$ be a simplicial sheaf with the properties that $X_m$ is a manifold for all $m\le n$,
and that it satisfies the Kan condition (\ref{dav}) for all $m\le n$.
Then if $B$ is a finite contractible simplicial set of dimension $\le n$, the sheaf $Hom(B,X)$ is representable by a manifold. \qed
\end{lemma}

\begin{proposition}\label{p6m}\dontshow{p6m}
Let $X$ be a reduced Kan simplicial manifold whose simplicial homotopy groups are finite dimensional diffeological groups.
Then $(\tau_{< n}X)_m$ is a manifold for every $m$, and the simplicial manifold $\tau_{< n}X$ satisfies the Kan condition.
\end{proposition}

\proof
For $m< n$, we have $(\tau_{< n}X)_m=X_m$ and so there's nothing to prove.
For $m=n$, the space $(\tau_{< n}X)_n$ is equal to $B_n$ and is a manifold by Lemma \ref{kz}.

Now we observe that in the diagram
\[
X_n\to B_n=(\tau_{< n}X)_n\to Hom(\Lambda[n,j],\tau_{< n} X)=Hom(\Lambda[n,j],X)
\]
the arrow 
$X_n\to B_n$ is surjective while the composite is a surjective submersion.
It follows from Lemma \ref{aA4} that
\[
B_n\to Hom(\Lambda[n,j],\tau_{< n} X)
\]
is a surjective submersion. 
We conclude that $\tau_{< n}X$ satisfies the Kan condition in dimensions $\le n$.
For $m>n$,
the space $(\tau_{< n}X)_m$ can be described explicitely as \dontshow{etl}
\begin{equation}\label{etl}
(\tau_{< n}X)_m = Hom(S,\tau_{< n}X),
\end{equation}
where $S\subset\Delta[m]$ is the union of all the
$n$-faces of $\Delta[m]$ incident to a given vertex;
it is then a manifold by Lemma \ref{KX}.

The Kan condition in dimensions $>n$ is trivial since the maps (\ref{dav}) are then always isomorphisms.
\qed

We can now finish the proof of our main Theorem.

\noindent {\it Proof of Theorem \ref{TVF}.}
For $m< n$, the sheaf $(\tau_{\le n}\intinf L)_m$ is equal to $(\intinf L)_m$, and is thus a manifold by Theorem \ref{tgo}.

Note that $(\tau_{\le n}\intinf L)_n$ is a $\pi_n^\simpl(\intinf L)$-principal bundle over $(\tau_{< n}\intinf L)_n$.
We know from Proposition \ref{p6m} that $(\tau_{< n}\intinf L)_n$ is a manifold.
Hence $(\tau_{\le n}\intinf L)_n$ is a manifold if and only if $\pi_n^\simpl(\intinf L)$ is a Lie group.
This is in turn the case if and only if $\partial_n:\pi_{n+1}(G)\to H_{n-1}(L)$ has discrete image.
Similarly, for $m>n$, the sheaf $(\tau_{\le n}\intinf L)_m$ is a  $\pi_n^\simpl(\intinf L)^{\left(\substack{\scriptstyle m\\ \scriptstyle n}\right)
}$-principal bundle over 
$(\tau_{< n}\intinf L)_m$ and is thus a manifold under the same assumption.

We now show that $\tau_{\le n} \intinf L$ satisfies the 
Kan condition.
For $m<n$, the map
\begin{equation}\label{ka}
(\tau_{\le n} \intinf L)_m\to Hom\big(\Lambda[m,j],\tau_{\le n} \intinf L\big)
\end{equation}
is a surjective submersion by Theorem \ref{tgo} since $\tau_{\le n} \intinf L=\intinf L$ in that range of dimensions.
For $m>n$, (\ref{ka}) is an isomorphism by Lemma \ref{check n-group}.
Finally, for $m=n$, the result follows by Theorem \ref{tgo} 
and Lemma \ref{aA4} since the composite
\[
(\intinf L)_n\twoheadrightarrow (\tau_{\le n} \intinf L)_n \to
Hom\big(\Lambda[n,j],\tau_{\le n} \intinf L\big)
=
Hom\big(\Lambda[n,j],\intinf L\big)
\]
is then a surjective submersion.
\qed

To complete the picture, we present an example where the hypotheses of Theorem \ref{TVF}
are not satisfied, and thus $\tau_{\le n}\intinf L$ is not composed of manifolds.

\begin{example}\label{eND}\dontshow{eND}
\rm
Let $\str$ be as in Definition \ref{d;}, let $p,q$ be $\mathbb Q$-linearly independent real numbers, and
let $L$ be the quotient of $\str\oplus\str$ given by $L_0=\g\oplus\g$ and $L_1=\R^2/(p,q)\R$.
By Lemma \ref{l,}, the boundary homomorphism $\partial:\pi_3(G\times G)\to H_1(L)$
is the composite of the standard inclusion $\Z^2\hookrightarrow \R^2$ with the projection
$\R^2\to\R^2/(p,q)\R$.
It doesn't have discrete image and so by Theorem \ref{TVF} the sheaf $(\tau_{\le 2}\intinf L)_2$ is not a manifold.

Indeed, anticipating (\ref{S1pb}), we can write $(\tau_{\le 2}\intinf (\str\oplus\str))_2$
as the total space of an $(S^1\times S^1)$-principal bundle.
The space 
$(\tau_{\le 2}\intinf L)_2$ is then the quotient of the above principal bundle by the non-closed subgroup
$(p,q)\R\subset S^1\times S^1$.
Note also that $\pi_2^\simpl(\intinf L)=\coker(\partial:\pi_3(G\times G)\to H_1(L))$ is not a Lie group.
\end{example}

\section{The string Lie 2-algebra}

This is our motivating example, which was introduced in \cite{BC04} and intensely studied in \cite{BCSS05}.

Let $\g$ be a simple real Lie algebra of compact type, 
and let $\langle\,,\rangle$ be the inner product on $\g$ such that the norm of its short coroots is 1.
\begin{definition}[\rm \cite{BCSS05}] \label{d;}\dontshow{d;}
Let $\g$ be a simple Lie algebra of compact type.
Its string Lie algebra is the Lie 2-algebra $\str=\str(\g)$ given by
\begin{equation*}\label{gg}
\str_0=\g, \qquad \str_1=\R
\end{equation*}
and brackets
\[
\begin{split}
[(X,c)]=0,\qquad [(X_1,c_1),(X_2,c_2)]=([X_1,X_2],0),\\
[(X_1,c_1),(X_2,c_2),(X_3,c_3)]=(0,\langle[X_1,X_2],X_3\rangle).\hspace{1.5mm}
\end{split}
\]
\end{definition}
The string Lie algebra should be thought as a central extension of the Lie algebra $\g$, 
but which is controlled by $H^3(\g,\R)$ as opposed to $H^2(\g,\R)$.
The Chevalley-Eilenberg complex of $\str$ is given by
\[
C^*(\str)=
\R\oplus \Big[\g^*\Big]\oplus \Big[\Lambda^2 \g^*\oplus \R\Big] \oplus \Big[\Lambda^3 \g^*\oplus \g^*\Big]\\
\oplus \Big[\Lambda^4 \g^*\oplus \Lambda^2 \g^*\oplus\R\Big]\oplus\ldots
\]
Following (\ref{in}), we study \dontshow{eqintstr}
\begin{equation}\label{eqintstr}
\begin{split}
Hom_{DGA}\big(C^*(\str),\Omega^*(\Delta^m)\big)=\big\{\alpha\in\Omega^1(\Delta^m;&\g),\beta\in\Omega^2(\Delta^m;\R)\,\big|\\
&d\alpha={\textstyle\frac{1}{2}}[\alpha,\alpha],\,
d\beta={\textstyle\frac{1}{6}}[\alpha,\alpha,\alpha]\big\}.\\
\end{split}
\end{equation}
The 1-form $\alpha$ satisfies the Maurer-Cartan equation, so we can integrate it to a map
$f:\Delta^m\to G$, defined up to translation, 
and satisfying $f^*(\theta_L)=-\alpha$, where $\theta_L\in\Omega^1(G;\g)$ is the left invariant Maurer-Cartan form on $G$.
The 3-form ${\textstyle\frac{1}{6}}[\alpha,\alpha,\alpha]$ is then the pullback of minus the Cartan 3-form 
\[
\eta=-{\textstyle\frac{1}{6}}\langle[\theta_L,\theta_L],\theta_L\rangle\in\Omega^3(G;\R),
\]
which represents a generator of $H^3(G,\Z)$.
So we can rewrite (\ref{eqintstr}) as \dontshow{es}
\begin{equation}\label{es}
\big(\intinf \str\big)_m=
\big\{f\in Map(\Delta^m,G)/G,\, \beta\in\Omega^2(\Delta^m)\,\big|d\beta=f^*(\eta)\big\}.
\end{equation}

In order to understand the simplicial homotopy groups of $\intinf \str$, 
we shall use the long exact sequence constructed in Theorem \ref{tls}.
First recall the well know facts $\pi_2(G)=0$, $\pi_3(G)=\Z$, which hold for any simple compact Lie group.
Let's also note that $H_0(\str)=\g$, $H_1(\str)=\R$ and $H_n(\str)=0$ for $n\ge 2$.
It follows that $\pi_n(\intinf\str)=\pi_n(G)$ for $n\ge 4$ and that we have an exact sequence
\dontshow{4Tm}
\begin{equation}\label{4Tm}
0\rightarrow\pi_3^\simpl(\intinf\str)\rightarrow\pi_3(G)\rightarrow H_1(\str)\rightarrow\pi_2^\simpl(\intinf\str)\rightarrow 0.
\end{equation}

\begin{lemma}\label{l,}\dontshow{l,}
The boundary homomorphism $\partial:\pi_3(G)\to H_1(\str)$ is the inclusion $\Z\hookrightarrow\R$.
\end{lemma}

\proof
Given an element $a\in\pi_3(G)=\Z$ represented by a map $f:\Delta^3\to G$ sending $\partial\Delta^3$ to the identity of $G$, we lift it to an element 
$(f,\beta)\in(\intinf \str)_3$.
When restricted to $\partial\Delta^3$, the map $f$ becomes constant and we're left with the 2-form $\beta|_{\partial\Delta^3}$.
This form then represents the element $\partial(a)=\int_{\partial\Delta^3}\beta\in\R=H_1(\str)$.
From Stokes' theorem, we then have
\[
\partial(a)=\int_{\partial\Delta^3}\beta=\int_{\Delta^3}f^*(\eta)=a
\]
as desired.
\qed

Thus we have proven:

\begin{theorem}
The simplicial homotopy groups of $\intinf\str$ are given by \dontshow{hkb}
\begin{equation}\label{hkb}
\begin{split}
\pi_1^\simpl(\intinf\str)=G,\quad\pi_2^\simpl(\intinf\str)&=S^1,\quad\pi_3^\simpl(\intinf\str)=0,\\ \pi_n^\simpl(\intinf\str)=\pi_n&(G)\quad \text{\rm for}\,\, n\ge 4.
\end{split}
\end{equation}
\vspace{-.4cm}

\qed
\end{theorem}

Moreover, it can be shown \cite{BMcL} that the first $k$-invariant of $\intinf\str$ is a generator of $H^3(K(G,1),S^1)=H^4(BG,\Z)=\Z$.
Here, $H^3(K(G,1),S^1)$ denotes simplicial manifold cohomology\footnote{Not to be confused with the continuous 
cohomology of $G$!} with coefficients in the sheaf of smooth $S^1$-valued functions 
(see \cite[Chapters 2,3]{Fri82} for definitions), 
while $H^4(BG,\Z)$ denotes usual topological cohomology of $BG=|K(G,1)|$.

We now study the Lie 2-group $\tau_{\le 2}\intinf\str$ integrating $\str$, and show that its geometric realization has the homotopy type of $B\Str$.
From (\ref{hkb}), we have
\begin{equation*}
\pi_1^\simpl(\tau_{\le 2}\intinf\str)=G\qquad\pi_2^\simpl(\tau_{\le 2}\intinf\str)=S^1,
\end{equation*}
and all other simplicial homotopy groups vanish.
The Postnikov tower (\ref{pKr}) thus induces a fiber sequence \dontshow{fS}
\begin{equation}\label{fS}
K(S^1,2)\,\rightarrow\,\tau_{\le 2}\intinf\str\,\rightarrow\,\tau_{< 2}\intinf\g.
\end{equation}
The base $\tau_{< 2}\intinf\str=\tau_{< 2}\intinf\g$ is
given by $(\tau_{< 2}\intinf\g)_m=Map(\sk_1 \Delta^m,G)/G$, and
the fiber is given by $\big(K(S^1,2)\big)_m=Z^2(\Delta[m],S^1)$, where $Z^2(\Delta[m],S^1)$ denotes the space of simplicial 2-cocycles on $\Delta[m]$ with values in $S^1$:
\[
K(S^1,2)=\left(*\llarrow *\lllarrow S^1\llllarrow (S^1)^3\lllllarrow (S^1)^6\cdots\right).
\]
Note that the fiber of the map $\tau_{< 2}\intinf\g\to \tau_{\le 1}\intinf\g=K(G,1)$ is simplicially contractible, and that therefore $|\tau_{< 2}\intinf\g\big|$ is a model for $BG$.
Also, by applying twice the simplicial path-loop fibration, 
one easily checks thats the geometric realization of $K(S^1,2)$ is a model for $K(\Z,3)$.
Upon geometric realization and up to homotopy, the fiber sequence (\ref{fS}) thus becomes \dontshow{KmK}
\begin{equation}\label{KmK}
K(\Z,3)\,\rightarrow\,\big|\tau_{\le 2}\intinf\str\big|\,\rightarrow\,BG.
\end{equation}
A version of the following Theorem has already appeared in \cite{BCSS05}. 
We include it here for the sake of completeness.

\begin{theorem}
The space $|\tau_{\le 2}\intinf\str\big|$ has the homotopy type of $B\Str$, where $\Str=\Str_G$ denotes the $3$-connected cover of $G$. 
\end{theorem}

\proof
Let $\tau_{>2}\intinf\str$ denote the fiber of the projection $\intinf\str\to\tau_{\le2}\intinf\str$.
Its simplicial homotopy groups are then given by \dontshow{4z}
\begin{equation}\label{4z}
\pi_n^\simpl(\tau_{>2}\intinf\str)=\begin{cases}
\pi_n^\simpl(\intinf\str)\quad &\text{for}\, n>2\\
0 &\text{otherwise.}
\end{cases}
\end{equation}
Since $\pi_n^\simpl(\tau_{>2}\intinf\str)$ 
are discrete groups, the simplicial Postnikov tower of $\tau_{>2}\intinf\str$ induces the usual Postnikov tower of $|\tau_{>2}\intinf\str|$
upon geometric realization.
It then follows from (\ref{hkb}) and (\ref{4z}) that $\pi_n|\tau_{>2}\intinf\str|=0$ for $n\le 3$.
The total space in the fiber sequence
\[
|\tau_{>2}\intinf\str| \to |\intinf\str| \to |\tau_{\le2}\intinf\str|
\]
is contractible as observed in section \ref{SEC6}. 
It follows that $\pi_n|\tau_{\le2}\intinf\str|=\pi_{n-1}|\tau_{>2}\intinf\str|$, and in particular that $\pi_n|\tau_{\le 2}\intinf\str|=0$ for $n\le 4$.

Since $\pi_4|\tau_{\le 2}\intinf\str|=0$, the fiber sequence (\ref{KmK}) identifies $|\tau_{\le 2}\intinf\str|$ with the 4-connected cover of $BG$. 
\qed

Our next task is to understand in explicit terms the manifolds composing $\inttwo \str$.
Note that for $m\le 1$, the manifold $(\inttwo\str)_m=(\intinf\str)_m$ has been computed in (\ref{es}).
So we concentrate on $(\inttwo\str)_2$.
Two 2-simplices $(f_0,\beta_0), (f_1,\beta_1)\in (\intinf\str)_2$ are identified in $(\inttwo\str)_2$ if they are simplicially homotopic relatively to $\partial\Delta[2]$.
This happens if 
there's a homotopy $F:\Delta^2\times[0,1]\to G$ from $f_0$ to $f_1$, fixing $\partial\Delta^2$,
and a 2-form $B$ on $\Delta^2\times[0,1]$ 
satisfying $B|_{\Delta^2\times\partial[0,1]}=\beta_0\sqcup\beta_1$ and
$dB=F^*\eta$.
By Stokes theorem, such a 2-form $B$ exists if and only if 
\[
\int_{\Delta^2\times[0,1]} F^*\eta=\int_{\Delta^2}\beta_1-\int_{\Delta^2}\beta_0.
\]
So we can forget about the actual 2-form $\beta$, and only keep its integral $b:=\int_{\Delta^2} \beta$.
We then get
\[
(\inttwo\str)_2=
\big(Map(\Delta^2,G)/G\big)
\times\R\Big/\sim\,,
\]
where $(f_0,b_0)$ and $(f_1,b_1)$ are identified if there exists a homotopy $F$ between $f_0$ and $f_1$, 
fixing $\partial\Delta^2$, and such that $\int F^*\eta=b_1-b_0$.
Note that since $F$ can wrap around $\pi_3(G)=\Z$, we have $(f,b)\sim(f,b+1)$,
so $(\inttwo\str)_2$ is the total space of an $S^1$-principal bundle \dontshow{S1pb}
\begin{equation}\label{S1pb}
S^1\,\to\,(\inttwo\str)_2\,\to\, Map(\partial\Delta^2,G)/G.
\end{equation}
The fiber of the map
\(
\big(Map(\Delta^2,G)/G\big)\times\R\rightarrow (\inttwo\str)_2
\)
is the universal cover of $\Omega^2 G= Map(\Delta^2/\partial\Delta^2,G)/G$. 
Since $\big(Map(\Delta^2,G)/G\big)\times\R$ is contractible, 
it follows that $(\inttwo\str)_2$ is homotopy equivalent to the 2-connected cover of $\Omega G$, and that the Chern class of (\ref{S1pb}) is a generator of 
$H^2(Map(\partial\Delta^2,G)/G,\Z)=H^2(\Omega G,\Z)=\Z$.

We now describe $(\inttwo\str)_m$ for arbitrary $m$.
Given a 2-face $\sigma:\Delta^2\to\Delta^m$, let $\cL_\sigma$ denote the pullback of (\ref{S1pb})
along the map $\sigma^*:Map(\sk_1\Delta^m,G)/G\to Map(\sk_1\Delta^2,G)/G$.
Given a 3-face $\tau:\Delta^3\to\Delta^m$, 
there is a preferred section $\phi_\tau$ of the tensor product
\begin{equation}\label{fgl}
\cL_{d_0\tau}\otimes\cL_{d_1\tau}^*\otimes\cL_{d_2\tau}\otimes\cL_{d_3\tau}^*,
\end{equation}
which we now describe.
Let us identify the dual of an $S^1$-bundle with the bundle itself equipped with the opposite $S^1$-action,
and let us call $a^*\in\cL^*$ the element corresponding to $a\in\cL$.
The section of (\ref{fgl}) is given by assigning to $f:\sk_1\Delta^m\to G$ the element
\begin{equation*}\label{fdx}
\phi_\tau(f):=
\big(d_0(
\tilde f\!\circ\!\tau
),b_0\big)\otimes\big(d_1(
\tilde f\!\circ\!\tau
),b_1\big)^*\otimes\big(d_2(
\tilde f\!\circ\!\tau
),b_2\big)\otimes\big(d_3(
\tilde f\!\circ\!\tau
),b_3\big)^*
\end{equation*}
defined as follows. 
Pick an extension $\tilde f:\Delta^m\to G$ of $f$ and let
$d_i(\tilde f \circ \tau):\Delta^2\to G$ be the $i$th face of $\tilde f \circ \tau$.
Let then $b_i\in\R$ be any numbers satisfying $b_0-b_1+b_2-b_3=\int_{\Delta^3}(\tilde f\circ\tau)^*\eta$.
The sections $\phi_\tau$ are well defined and satisfy the cocycle condition saying that, for every 4-face $\nu:\Delta^4\to\Delta^m$ the section
$\phi_{d_0\nu}\otimes\phi_{d_1\nu}^*\otimes\phi_{d_2\nu}\otimes\phi_{d_3\nu}^*\otimes\phi_{d_4\nu}$
is equal to the trivial section in the appropriate trivialized $S^1$-bundle.

Given the above structure on the $S^1$-bundles $\cL_\sigma$, 
we now have the following description of $(\inttwo\str)_m$.

\begin{proposition}
Let $\cL_\sigma$ be the above $S^1$-bundle over $Map(\sk_1(\Delta^m,G)/G)$, 
and let $\cL_{\sigma,f}$ denote its fiber over a given point $f$. 
We then have \dontshow{frzx}
\begin{equation}\label{frzx}
(\inttwo\str)_m=\Big\{
f\in Map(\sk_1\Delta^m\!,G)/G,\,\,(c_\sigma\in\cL_{\sigma,f})\,\Big|\,\delta c=\phi
\Big\},
\end{equation}
where $\delta c=\phi$ means that for every 3-face $\tau:\Delta^3\to\Delta^m$ we have \dontshow{dzx}
\begin{equation}\label{dzx}
c_{d_0\tau}\otimes c_{d_1\tau}^*\otimes c_{d_2\tau}\otimes c_{d_3\tau}^*=\,\phi_\tau\,.
\end{equation}
\end{proposition}

\proof
Recall that $(\inttwo\str)_m$ is a quotient of $(\intinf\str)_m$ where two elements
are identified if they are simplicially homotopic relatively to $\sk_1\Delta[m]$.
As explained in (\ref{es}), an element of $(\intinf\str)_m$ is represented by a pair $(f,\beta)$ where $f$ is a map $\Delta^m \to G$ and $\beta$ is a 2-form 
on $\Delta^m$ satisfying $d\beta=f^*\eta$.
Integrating $\beta$ over the various 2-faces $\sigma$ of $\Delta^m$ produces real numbers
$b_\sigma\in \R$, and the equivalence class of $(f,\beta)$ only depends on $(f,(b_\sigma))$.
By Stokes' theorem, we then have for every 3-face $\tau$ the relation \dontshow{dzl}
\begin{equation}\label{dzl}
b_{d_0\tau}-b_{d_1\tau}+b_{d_2\tau}-b_{d_3\tau}=\int_{\Delta^3}\tau^*f^*\eta.
\end{equation}

Given $(f,(b_\sigma))$ representing an element in the left hand side of (\ref{frzx}), 
we can assign to it the pair $(f|_{\sk_1\Delta^m}, (c_{\sigma}))$,
where $c_\sigma\in\cL_{\sigma,f}$ is the element represented by
$(d_i(f \circ\tau),b_i)$.
Equation (\ref{dzl}) then implies (\ref{dzx}) by definition of $\phi_\tau$.

Inversely, if we are given $(f,(c_\sigma))$ in the right hand side of (\ref{frzx}), we first pick an extension $\tilde f:\Delta^m\to G$.
Let $S\subset \Delta^m$ be the union
of all 2-faces incident to a given vertex. 
It is a maximal contractible subcomplex of $\sk_2\Delta^m$.
Since $S$ is contractible, there are no obstructions to picking a 2-form on $\Delta^m$ with prescribed integrals on the various 2-faces of $S$.
So we may pick $\beta\in\Omega^2(\Delta^m)$ such that for every face 
$\sigma:\Delta^2\to S$, the pair $\big(\tilde f\circ \sigma,\,\int_{\Delta^2}\sigma^*\beta\big)$ is a representative of $c_\sigma$.
From (\ref{dzx}) and Stokes' theorem, we conclude that $\big(\tilde f\circ \sigma,\,\int_{\Delta^2}\sigma^*\beta\big)$
is also a representative of $c_\sigma$ for the other 2-faces $\sigma$ of $\Delta^m$.
It follows that the assignment $(f,(c_\sigma))\mapsto (\tilde f,\beta)$ is an inverse to the previous construction.
\qed

\begin{corollary}
The manifold $(\inttwo\str)_m$ is the total space of an $(S^1)^{\left(\substack{\scriptstyle m\\ \scriptstyle 2}\right)
}$-principal bundle over $Map(\sk_1\Delta^m,G)/G$.
\end{corollary}

Note that the abelian group $(S^1)^{\left(\substack{\scriptstyle m\\ \scriptstyle 2}\right)}$ can be naturally identified with 
$S^1\otimes H^1(sk_1 \Delta^m,\Z)$.
So we can write schematically our computation of $\inttwo\str$ as \dontshow{TheBstr}
\begin{equation}\label{TheBstr}
\inttwo\str=\left[*\llarrow Path(G)/G\lllarrow \widetilde{Map(\partial \Delta^2, G)}/G
\llllarrow \widetilde{Map(sk_1 \Delta^3, G)}/G \cdots\right]
\end{equation}
where the tilde indicates that $Map(sk_1 \Delta^m, G)$ has been replaced by the total space of an 
$(S^1\otimes H^1(sk_1 \Delta^m,\Z))$-principal bundle.

If we apply the simplicial loop functor to (\ref{TheBstr}), we get another simplicial manifold whose homotopy type is now that of $\Str_G$:
\[\Omega^\simpl(\inttwo\str)=\left[
Path(G)/G\llarrow 
\widetilde{Loop(G)}/G\lllarrow
\widetilde{Map(sk_1(\Sigma\Delta^2),G)}/G\llllarrow\cdots\right],
\]
or more pictorially \dontshow{TheStr}
\begin{equation}\label{TheStr}
Path_*(G)\llarrow 
\widetilde{Loop_*(G)} \lllarrow
\widetilde{Map_*(\put(5,3){\circle{8}}\put(1,3){\line(1,0){8}}\hspace{.4cm},G)}
\llllarrow
\widetilde{Map_*(
\put(5.1,3){\circle{8.5}}\put(.2,-.5){$\smile$}\put(.2,1.3){$\frown$}\hspace{.4cm}
,G)}
\cdots,
\end{equation}
where the stars denote pointed mapping spaces.

\begin{remark} \rm 
It was a pleasant surprise to see that the manifolds in (\ref{TheStr}) come with natural group structures, 
and that these group structures are preserved by the face and degeneracy maps.
Therefore (\ref{TheStr}) is a simplicial Banach Lie group, and its geometrical realization is a topological group.
We don't know if that's a general phenomenon, or if it is a special feature of the Lie 2-algebra $\str$.

The manifolds in (\ref{TheBstr}) can also be endowed with group structures, but not all face maps are group homomorphisms.
\end{remark}

\section{Appendix}

In this appendix, we compare our notion of Lie 2-groups with the notion of coherent Lie 2-group introduced by Baez and Lauda \cite[Definition 19]{BL04}.
To avoid confusions, we call their notion a Lie $2_{B\!L}$-group and ours a Lie $2_H$-group.

A Lie $2_{B\!L}$-group is similar to a group object in the category of Lie groupoids.
It consists of a Lie groupoid $G=(G_1\rrarrow G_0)$, a multiplication
$(\mu:G^2\to G)=(\mu_0:G_0^2\to G_0;\mu_1:G_1^2\to G_1)$, and an associator 
$\alpha:G_0^3\to G_1$, which is a natural transformation between 
$\mu\circ(1\times \mu)$ and $\mu\circ(\mu\times 1)$, and satisfies the pentagon identity. 
Some other pieces of data describe the unit and inverse maps.
The Lie $2_H$-group corresponding to $G$ is then given by
\[
\begin{split}
NG_m:=\Big\{(g_{ij}\in G_0)_{0\le i<j\le m},\,(h_{ijk}:g_{ij}g_{jk}\to g_{ik})_{0\le i<j<k\le m}\,\Big|\,
&\forall i\!<\!j\!<\!k\!<\!\ell,\\ &\text{(\ref{ctt}) commutes}
\Big\}
\end{split}
\]
\begin{equation}\label{ctt}
\begin{matrix}\xymatrix{
g_{ij}(g_{jk}g_{k\ell})\ar[d]_{g_{ij}h_{jk\ell}}\ar[rr]^{\alpha(g_{ij},g_{jk},g_{k\ell})}&&(g_{ij}g_{jk})g_{k\ell}
\ar[d]^{h_{ijk}g_{k\ell}}\\
g_{ij}g_{j\ell}\ar[r]^{h_{ij\ell}}&g_{i\ell}&\ar[l]_{h_{ik\ell}}g_{ik}g_{k\ell}
}\end{matrix}
\end{equation}
The only non-trivial thing to check is that 
any map $\Lambda[3,j]\to NG$ extends uniquely to a map $\Delta[3]\to NG$.
This follows from the fact that any one of the $h$'s in (\ref{ctt}) is uniquely determined by the three other ones.

For the inverse construction, we first do the discrete case.
Given a discrete $2_H$-group $X$, we consider its simplicial loop space $\Omega X$
\[
(\Omega X)_m=\{x\in X_{m+1}\,|\,d_{m+1}(x)=*\}.
\]
The $2_{B\!L}$-group $G$ associated to $X$ is then the fundamental groupoid of $\Omega X$
\begin{equation}\label{DEFG}
G_0=X_1,\qquad G_1=\{x\in X_2\,|\,d_2(x)=*\}.
\end{equation}
The multiplication on objects $g,g'\in G_0$ is given by picking a filler $x=x(g,g') \in X_2$ such that $d_0(x)=g$, $d_2(x)=g'$ and letting $gg':=d_1(x)$. 
The other pieces of data are then constructed by similar but more complicated horn filling procedures.

For Lie $2_H$-groups, the above construction is unfortunately problematic.
We can define $G$ as in (\ref{DEFG}), but the multiplication $\mu_0:G_0\times G_0\to G_0$ requires
the choice of a filler $x(g,g')$, and that choice can only be made {\em locally} in $(g,g')$.
Our opinion is that the notion of coherent Lie 2-group presented in \cite{BL04} is maybe ``not weak enough''.
A weaker notion might look as follows.

Given a Lie groupoid $G$ and a cover $U\to G_0$ of its manifold of objects, 
let $G_U$ denote the groupoid $U\times_{G_0}G_1\times_{G_0}U\rrarrow U$.

\begin{definition}\label{soD}
A weak Lie 2$_{B\!L}$-group is a Lie groupoid $G=(G_1\rrarrow G_0)$
equipped with the following pieces of data:

$\bullet$ An identity element $e\in G_0$.

$\bullet$ A cover $p:U\to G^2_0$, and two multiplications maps
$$\mu_0:U\to G_0,\qquad \mu_1:U\times_{G^2_0}G^2_1\times_{G^2_0}U\to G_1$$ 
assembling to a smooth functor $\mu$ from $(G^2)_U$ to $G$.

$\bullet$ Two smooth maps 
\[
\begin{split}
\ell:\,_eU&=\big\{(g,u)\in G_0\times U\,|\, p(u)=(e,g)\big\}\to G_1,\\
r:U_e&=\big\{(g,u)\in G_0\times U\,|\, p(u)=(g,e)\big\}\to G_1,
\end{split}
\]
forming natural transformations
\[
\xymatrix{
&(G^2)_U\ar[dr]^\mu&\\
\,\,G_{_eU}\ar[ur]^{e\times 1}\ar[rr]&\ar@{}[u]|(.38){\,\,\,\displaystyle \Downarrow \ell}&G\,,
}
\hspace{1cm}
\xymatrix{
&(G^2)_U\ar[dr]^\mu&\\
\,\,G_{U_e}\ar[ur]^{1\times e}\ar[rr]&\ar@{}[u]|(.38){\,\,\,\displaystyle \Downarrow r}&G\,.
}
\]

$\bullet$ An associator
\[
\begin{split}
\alpha:V\!=\!\big\{(g_0&,g_1,g_2,u_{01},u_{12},u',u'')\in G_0^3\times U^4\,\big|\,p(u_{ij})=(g_i,g_j),\\
&p(u')=(g_0,\mu(u_{12})),\,p(u'')=(\mu(u_{01}),g_3)\big\}\to G_1
\end{split}
\]
which is a natural transformation 
\[
\xymatrix@R=.4cm{
&
\big(G\times(G^2)_U\big)_{U'}\ar[r]^(.6){1\times \mu}& (G^2)_U\ar[dr]^(.55)\mu&\\
\save[]+<.3cm,0cm>*\txt<2pc>{$(G^3)_V$\phantom{\Big|}}\ar[ur]\ar[dr]\restore
\ar@{}[rrr]|(.55){\displaystyle\Downarrow\alpha}&&&G\,\,,\\
&\big((G^2)_U\times G\big)_{U''}\ar[r]^(.6){\mu\times 1}& (G^2)_U\ar[ur]_(.55)\mu&
}
\]
where $U'$ and $U''$ denote the covers $(1\times \mu_0)^{-1}(U)$ and $(\mu_0\times 1)^{-1}(U)$ respectively.
\medskip

\noindent
The maps $e,\mu_0,\mu_1,\ell, r,\alpha$ are subject to the following three axioms:
\medskip

$\bullet$ 
Letting $W$ be the cover of $G_0^4$ given by
\[
\begin{split}
W\!=\!\big\{
(g_0,g_1,g_2,g_3,&u_{01},u_{12},u_{23},u'_{012},u''_{012},u'_{123},u''_{123},v_1,v_2,v_3,v_4,v_5)\in G_0^4\times U^{12}\,\big|\\
p(u_{ij})\!&=\!(g_i,g_j),\,
p(u'_{ijk})\!=\!(g_i,\mu(u_{jk})),\,
p(u''_{ijk})\!=\!(\mu(u_{ij}),g_k),\\
p(v_1)\!&=\!(g_0,\mu(u'_{123})),\,
p(v_2)\!=\!(\mu(u_{01}),\mu(u_{23})),\\
p(v_3)\!&=\!(\mu(u''_{012}),g_3),\,
p(v_4)\!=\!(\mu(u'_{012}),g_3),\,
p(v_5)\!=\!(g_0,\mu(u''_{123}))
\big\},
\end{split}
\]
then the five
associators between the five possible multiplications $(G^4)_W\to G$ form a commutative diagram of natural transformations.

$\bullet$ Letting $V_e$ be the cover of $G_0^2$ given by
\[
\begin{split}
V_e=\big\{(g_0,g_1,u_0,u_1,u'_0,u'_1,u'')\in G_0^2\times U^4\,\big|\, p&(u_0)=(g_0,e),\,p(u_1)=(e,g_1),\\p(u_0')=(\mu(u_0),g_1),\,&p(u_1')=(g_0,\mu(u_1)),\,p(u'')=(g_0,g_1)
\big\},
\end{split}
\]
then the three natural transformations induced by $\alpha$, $\ell$, $r$, 
between the three possible multiplications $(G^2)_{V_e}\to G$ form a commutative diagram.
\medskip

$\bullet$ 
For each point $g\in G_0$, there exists a neighborhood $V(g)\subset G_0$, and smooth maps 
\[
\begin{split}
i_g&:V(g)\to G_0,\\ 
\lambda_g&:\{(x,u)\in V(g)\times U\,|\,p(u)=(x,i_g(x))\}\to G_1,\\
\rho_g&:\{(x,u)\in V(g)\times U\,|\,p(u)=(i_g(x),x)\}\to G_1,
\end{split}
\]
such that $\lambda_g(x,u)$ and $\rho_g(x,u)$ are arrows with source $\mu(u)$ and target $e$.
\end{definition}

\begin{remark}\rm
If we're only interested about what happens in a neighborhood of the identity, then Definition \ref{soD} is not really needed 
since the choice of fillers can be made locally.
It follows that germs of Lie $2_{B\!L}$-group are equivalent to germs of Lie $2_H$-groups.
\end{remark}

\bibliography{../main}
\bibliographystyle{plain}

\end{document}